\documentclass[a4paper,12pt]{article}
\usepackage{amsmath}
\usepackage{amsfonts}
\usepackage{amssymb}
\usepackage{latexsym}
\usepackage{epsfig}
\usepackage{graphicx}
\usepackage{oldgerm}
\usepackage{theorem}

\def\C{\mathbb{C}}
\def\R{\mathbb{R}}
\def\Z{\mathbb{Z}}
\def\N{\mathbb{N}}

\def\P{{\mathbb P}}

\def\E{{\mathbb E}}

\def\mbK{\mathbb{K}}

\def\T{{\bf T}}
\def\bK{{\bf K}}
\def\bP{{\bf P}}

\def\X{\mib{X}}
\def\Y{\mib{Y}}
\def\x{\mib{x}}

\def\0{\mib{0}}

\def\mM{\mathfrak{M}}
\def\mX{\mathfrak{X}}
\def\mY{\mathfrak{Y}}

\def\supp{{\rm supp}\ }

\def\Ai{{\rm Ai}}
\def\cA{{\cal A}}
\def\cK{{\cal K}}
\def\cG{{\cal G}}
\def\cS{{\cal S}}

\theorembodyfont{\itshape}

\newtheorem{thm}{Theorem}[section]
\newtheorem{lem}[thm]{Lemma}

\newtheorem{prop}[thm]{Proposition}

\newcommand{\mib}[1]{\mbox{\boldmath $#1$}}
\newcommand{\SSC}[1]{\section{#1}\setcounter{equation}{0}}
\newcommand{\qed}{\hbox{\rule[-2pt]{3pt}{6pt}}}

\begin{document}

\title{Zeros of Airy Function and Relaxation Process}
\author{
Makoto Katori
\footnote{
Department of Physics,
Faculty of Science and Engineering,
Chuo University, 
Kasuga, Bunkyo-ku, Tokyo 112-8551, Japan;
e-mail: katori@phys.chuo-u.ac.jp
}
and 
Hideki Tanemura
\footnote{
Department of Mathematics and Informatics,
Faculty of Science, Chiba University, 
1-33 Yayoi-cho, Inage-ku, Chiba 263-8522, Japan;
e-mail: tanemura@math.s.chiba-u.ac.jp
}}
\date{29 September 2009}
\pagestyle{plain}
\maketitle
\begin{abstract}
One-dimensional system of Brownian motions called
Dyson's model is the particle system with long-range repulsive
forces acting between any pair of particles,
where the strength of force is
$\beta/2$ times the inverse of particle distance.
When $\beta=2$, it is realized as the Brownian motions
in one dimension conditioned never to collide with
each other.
For any initial configuration, it is proved that
Dyson's model with $\beta=2$ and $N$ particles,
$\X(t)=(X_1(t), \dots, X_N(t)), 
t \in [0,\infty), 2 \leq N < \infty$,
is determinantal in the sense that any multitime correlation
function is given by a determinant with a continuous kernel.
The Airy function $\Ai(z)$ is an entire function
with zeros all located on the negative part of the 
real axis $\R$. We consider Dyson's model with $\beta=2$ 
starting from the first $N$ zeros of $\Ai(z)$,
$0 > a_1 > \cdots > a_N$, $N \geq 2$.
In order to properly control the effect of such 
initial confinement of particles in the negative region
of $\R$,
we put the drift term to each Brownian motion,
which increases in time as a parabolic function : 
$Y_j(t)=X_j(t)+t^2/4+\{d_1+\sum_{\ell=1}^{N}(1/a_{\ell})\}t,
1 \leq j \leq N$,
where $d_1=\Ai'(0)/\Ai(0)$.
We show that, as the $N \to \infty$ limit of
$\Y(t)=(Y_1(t), \dots, Y_N(t)), t \in [0, \infty)$, we obtain
an infinite particle system, which is the relaxation
process from the configuration, in which every zero of $\Ai(z)$
on the negative $\R$ is occupied by one particle, 
to the stationary state $\mu_{\Ai}$.
The stationary state $\mu_{\Ai}$ is the determinantal point
process with the Airy kernel, which is spatially inhomogeneous
on $\R$ and in which the Tracy-Widom distribution 
describes the rightmost particle position.

\end{abstract}

\clearpage

\SSC{Introduction}\label{chap: Introduction}

\subsection{Dyson's model: One-Dimensional
Brownian Particle System Interacting through
Pair Force $1/x$}

To understand the time-evolution of distributions of
interacting particle systems on a large space-time
scale ({\it thermodynamic and hydrodynamic limits})
is one of the main topics of statistical physics.
If the interactions among particles are short ranged,
the standard theory is useful.
If they are long ranged, however,
general theory has not yet been established
and thus detailed study of model systems 
is required \cite{Spo91}.

In the present paper, we consider Brownian particles
in one dimension with long-ranged repulsive forces
acting between any pair of particles,
where the strength of force is exactly equal to
$1/x$ when the particle distance is $x$.
If the number of particles is finite $N < \infty$,
the system is described by
$\Xi(t)=\sum_{j=1}^{N} \delta_{X_j(t)}$, 
$2 \leq N < \infty$, 
where $\X(t)=(X_1(t), \dots, X_N(t))$ satisfies the
following system of stochastic differential equations (SDEs);
\begin{equation}
dX_j(t)=dB_j(t) 
+ 
\sum_{\substack{1 \leq k \leq N
\\ k \not= j}}
\frac{1}{X_j(t)-X_k(t)} dt, \quad
1 \leq j \leq N, \quad
t \in [0, \infty)
\label{eqn:Dyson}
\end{equation}
with independent one-dimensional
standard Brownian motions $B_j(t), 1 \leq j \leq N$.
The SDEs obtained by replacing the $1/x$ force 
in (\ref{eqn:Dyson}) by $\beta/(2x)$ with
a parameter $\beta >0$ were 
introduced by Dyson \cite{Dys62}
to understand the statistics of eigenvalues of hermitian
{\it random matrices} as particle distributions of interacting 
Brownian motions in $\R$.
Corresponding to the special values
$\beta=1,2$ and 4, hermitian random matrices are
in the three statistical ensembles with different symmetries, 
called the Gaussian orthogonal ensemble (GOE), 
the Gaussian unitary ensemble (GUE),
and the Gaussian symplectic ensemble (GSE),
respectively \cite{Dys62b}.
In particular for $\beta=2$, 
that is the case of (\ref{eqn:Dyson}), 
if the eigenvalue distribution of
$N \times N$ hermitian random matrices in the GUE
with variance $\sigma^2$ 
is denoted by $\mu^{\rm GUE}_{N, \sigma^2}$,
we can show
\begin{equation}
\lim_{N \to \infty} \mu^{\rm GUE}_{N, 2N/\pi^2} 
(\cdot)
= \mu_{\sin}(\cdot),
\label{eqn:lim_sin}
\end{equation}
where $\mu_{\sin}$ denotes the 
{\it determinantal (Fermion) point process} \cite{Sos00,ST03}
with the so-called {\it sine kernel}
\begin{equation}
K_{\sin}(x)= 
\frac{1}{2 \pi} \int_{|k| \leq \pi} dk \,
e^{\sqrt{-1} k x}
= \frac{\sin (\pi x) }{\pi x},
\quad x \in \R.
\label{eqn:Ksin}
\end{equation}
That is, $\mu_{\sin}$ is a spatially homogeneous 
particle distribution, in which
the particle density is given by
$\rho_{\sin}=\lim_{x \to 0} K_{\sin}(x)=1$ and 
any $N_1$-point correlation function
$\rho_{\sin}(\x_{N_1}), \x_{N_1}=(x_1, \dots, x_{N_1}) \in \R^{N_1}, 
N_1 \geq 2$,
is given by a determinant of an $N_1 \times N_1$
real symmetric matrix;
$$
\rho_{\sin}(\x_{N_1})= \det_{1 \leq j, k \leq N_1} 
\Big[ K_{\sin}(x_j-x_k) \Big].
$$
Based on this fact known in the random matrix theory \cite{Meh04},
Spohn \cite{Spo87} studied the equilibrium dynamics 
obtained in the infinite-particle limit $N \to \infty$ 
of Dyson's model (\ref{eqn:Dyson}).
Since the $1/x$ force is not summable,
in the infinite-particle limit $N \to \infty$
the sum in (\ref{eqn:Dyson})
should be regarded as an improper sum,
in the sense that for $X_j(t) \in [-L,L]$
the summation is restricted to $k$'s such that
$X_k(t) \in [-L, L]$ and then
the limit $L \to \infty$ is taken.
It is expected that the dynamics with an infinite number
of particles can exist only for
initial configurations having the same
asymptotic density to the right and left
\cite{Spo87,KT08}.

The problem, which we address in the present paper,
is how we can control Dyson's model with 
an infinite number of particles starting from
{\it asymmetric} initial configurations.
The motivation is again coming from 
the random matrix theory as follows.
Consider the {\it Airy function} \cite{AS65,VS04}
\begin{equation}
\Ai(z) = \frac{1}{2 \pi} \int_{\R} dk \,
e^{\sqrt{-1}(z k+k^3/3)}.
\label{eqn:Airy1}
\end{equation}
It is a solution of Airy's equation
$f''(z)-z f(z)=0$
with the asymptotics on the real axis $\R$:
\begin{eqnarray}
&& \Ai(x) \simeq \frac{1}{2 \sqrt{\pi} x^{1/4}}
\exp\left( - \frac{2}{3} x^{3/2} \right), 
\nonumber\\
&& \Ai(-x) \simeq \frac{1}{\sqrt{\pi} x^{1/4}}
\cos \left( \frac{2}{3}x^{3/2}-\frac{\pi}{4} \right)
\quad \mbox{in} \quad x \to + \infty.
\label{eqn:Aiasym}
\end{eqnarray}
In the GUE random matrix theory,
the following scaling limit 
has been extensively studied:
\begin{equation}
\lim_{N \to \infty} \mu^{\rm GUE}_{N, N^{1/3}}
(2N^{2/3}+\cdot)=\mu_{\Ai}(\cdot),
\label{eqn:lim_Ai}
\end{equation}
where $\mu_{\Ai}$ is the determinantal point process
such that the correlation kernel is given by \cite{For93,TW94},
\begin{eqnarray}
K_{\Ai}(y|x)
&=& \int_{0}^{\infty} du \,
\Ai(u+x) \Ai(u+y) \nonumber\\
\label{eqn:AiryKer}
&=& \left\{ \begin{array}{ll}
\displaystyle{
\frac{\Ai(x) \Ai'(y)-\Ai'(x) \Ai(y)}{x-y}},
& x \not= y \in \R \cr
(\Ai'(x))^2-x (\Ai(x))^2,
& x=y \in \R. \end{array} \right.
\end{eqnarray}
It is another infinite-particle limit different from 
(\ref{eqn:lim_sin}) and is called
the {\it soft-edge scaling limit},
since $x^2/2t \simeq (2 N^{2/3})^2/(2N^{1/3})=2N$
marks the right edge of semicircle-shaped
profile of the GUE eigenvalue distribution
(see, for example, \cite{KT07}).
The particle distribution $\mu_{\Ai}$ with
the {\it Airy kernel} (\ref{eqn:AiryKer}) is
highly asymmetric:
As a matter of fact, the particle
density $\rho_{\Ai}(x)=K_{\Ai}(x|x)$
decays rapidly to zero as $x \to \infty$, 
but it diverges
\begin{equation}
\rho_{\Ai}(x)
\simeq \frac{1}{\pi} (-x)^{1/2} \to \infty
\quad \mbox{as} \quad x \to - \infty.
\label{eqn:rhoAi}
\end{equation}
Let $R$ be the position of the rightmost particle 
on $\R$ in $\mu_{\Ai}$.
Then its distribution is given by the
celebrated {\it Tracy-Widom distribution} \cite{TW94}
$$
\mu_{\Ai}(R < x)=\exp \left[
-\int_{x}^{\infty}(y-x) (q(y))^2 dy \right],
$$
where $q(x)$ is the unique solution of the Painlev\'e II equation
$$
q''=xq + 2 q^3
$$
satisfying the boundary condition
$q(x) \simeq \Ai(x)$ in $x \to \infty$.
Pr\"ahofer and Spohn \cite{PS02} and Johansson \cite{Joh03}
studied the equilibrium fluctuation of this rightmost particle 
and called it the {\it Airy process}.
Tracy and Widom derived a system of
partial differential equations, 
which govern the Airy process \cite{TW03}.
See also \cite{AM99, AM05}.
How can we realize $\mu_{\Ai}$ as the equilibrium state
of Brownian infinite-particle system interacting
through pair force $1/x$ ?
The initial configurations should be asymmetric,
but what kinds of conditions should be satisfied
by them ?
How should we modify the SDEs of original Dyson's model
(\ref{eqn:Dyson}), when we provide finite-particle
approximations for such asymmetric infinite particle
systems ?

In the present paper, as an explicit 
answer to the above questions, 
we will present a relaxation process
with an infinite number of particles converging
to the stationary state $\mu_{\Ai}$ in $t \to \infty$.
Its initial configuration is given by
\begin{equation}
\xi_{\cA}(\cdot)=\sum_{a \in {\cA} } \delta_a (\cdot)
=\sum_{j=1}^{\infty} \delta_{a_{j}}(\cdot),
\label{eqn:xiA}
\end{equation}
in which every zero of
the Airy function (\ref{eqn:Airy1}) 
is occupied by one particle.
This special choice of the initial configuration is 
due to the fact that
the zeros of the Airy function 
are located only on 
the negative part of the real axis $\R$, 
\begin{equation}
\cA \equiv \Ai^{-1}(0)
=\Big\{ a_{j}, j \in \N \, : \, \Ai(a_{j})=0, \,
0 > a_1 > a_2 > \cdots \Big\},
\label{eqn:AiryZero}
\end{equation}
with the values \cite{AS65}
$
a_1 = -2.33 \dots, 
a_2 = -4.08 \dots, 
a_3 = -5.52 \dots, 
a_4 = -6.78 \dots,
$
and that they
admit the asymptotics \cite{AS65,VS04}
\begin{equation}
a_{j} \simeq - \left( \frac{3 \pi}{2} \right)^{2/3} j^{2/3} \quad
\mbox{in} \quad j \to \infty.
\label{eqn:ainf}
\end{equation}
Then the average density of zeros
of the Airy function around $x$,
denoted by $\rho_{\Ai^{-1}(0)}(x)$, behaves as
$$
\rho_{\Ai^{-1}(0)}(x) \simeq \frac{1}{\pi} (-x)^{1/2}
\to \infty
\quad \mbox{as} \quad x \to -\infty,
$$
which coincides with (\ref{eqn:rhoAi}).
The approximation of our process
with a finite number of particles $N < \infty$ is
given by 
$\Xi_{\cA}(t)=\sum_{j=1}^{N} \delta_{Y_j(t)}$ with
\begin{equation}
Y_j(t) = X_j(t)+\frac{t^2}{4} + D_{\cA_N} t, \quad
1 \leq j \leq N, \quad
t \in [0, \infty),
\label{eqn:finite_Airy}
\end{equation}
associated with the solution
$\X(t)=(X_1(t), \dots, X_N(t))$ of Dyson's model 
(\ref{eqn:Dyson}),
where
\begin{equation}
D_{\cA_N}=d_1+\sum_{\ell=1}^N \frac{1}{a_\ell}.
\label{eqn:DA}
\end{equation}
Here
$d_1=\Ai'(0)/\Ai(0)$ and
$
{\cA}_N \equiv \Big\{0 > a_1 > \dots > a_N \Big\} 
\subset {\cA}
$
is the sequence of the first $N$ zeros
of the Airy function.
In other words, $\Y(t)=(Y_1(t), Y_2(t), \dots, Y_N(t))$
satisfies the following SDEs ;
\begin{eqnarray}
dY_j(t) &=& dB_j(t)+
\left( \frac{t}{2}+D_{\cA_N} \right) dt
+\sum_{\substack{1 \leq k \leq N
\\ k \not= j}}
\frac{dt}{Y_j(t)-Y_k(t)}
\nonumber\\
&=& dB_j(t)+\sum_{\substack{1 \leq k \leq N
\\ k \not= j}}
\left( \frac{1}{Y_j(t)-Y_k(t)}+\frac{1}{a_k}\right)dt
+ \left(\frac{t}{2} +d_1+\frac{1}{a_j}\right)dt,
\nonumber\\
&& \hskip 5cm 1 \leq j \leq N, \quad
t \in [0, \infty),
\label{eqn:YN}
\end{eqnarray}
where $B_j(t)$'s are independent one-dimensional
standard Brownian motions.

For $\Y(0)=\x \in \R^N$, set
$\xi^{N}(\cdot)=\sum_{j=1}^{N} \delta_{x_j}(\cdot)$
and consider the process $\Xi_{\cA}(t)$ starting from
the configuration $\xi^N$.
We consider a set of initial configurations $\xi^N$
such that they are in general different from
the $N$-particle approximation
of (\ref{eqn:xiA}),
\begin{equation}
\xi_{\cA}^N(\cdot)
= \sum_{a \in {\cA}_N} \delta_a(\cdot)
=\sum_{j=1}^{N} \delta_{a_j}(\cdot),
\label{eqn:xiAN}
\end{equation}
but the particle density $\rho(x)$
of $\lim_{N \to \infty} \xi^N$ 
will show the same
asymptotic in $x \to -\infty$ as (\ref{eqn:rhoAi}).
Because of the strong repulsive forces acting
between particle pairs in (\ref{eqn:Dyson}),
such confinement of particles in the negative region
of $\R$ at the initial time causes 
strong {\it positive drifts} of Brownian particles.
The coefficient (\ref{eqn:DA}) of the drift term 
$D_{\cA_N} t$ added in
(\ref{eqn:finite_Airy}), however, {\it negatively diverges}
\begin{equation}
D_{{\cA}_N} \simeq - \left(\frac{12}{\pi^2} \right)^{1/3}
N^{1/3} \to - \infty
\quad \mbox{as} \quad N \to \infty.
\label{eqn:Dinf}
\end{equation}
We will determine a class of asymmetric initial
configurations denoted by $\mX^{\cA}_0$, 
which includes $\xi_{\cA}$
as a typical one, such that the effect on dynamics 
of asymmetry in configuration will be compensated by
the additional drift term $D_{\cA_N} t$ in the
infinite-particle limit $N \to \infty$
and the dynamics with an infinite number of particles
exists.
Note that we should take $N \to \infty$ limit for 
finite $t < \infty$ in our process (\ref{eqn:finite_Airy})
to discuss non-equilibrium dynamics 
with an infinite number of particles. 
In the class $\mX^{\cA}_0$, when 
the initial configuration is specially
set to be $\xi_{\cA}$, (\ref{eqn:xiA}), 
we can prove that 
the dynamics shows a relaxation 
in the long-term limit $t \to \infty$
to the equilibrium dynamics in $\mu_{\Ai}$.
In the proof we use the special property of the systems
(\ref{eqn:Dyson}) and (\ref{eqn:YN}) such that
the processes have space-time determinantal
correlations.
This feature comes from the fact that
if and only if the strength of pair force
is exactly equal to $1/x$ when the particle
distance is $x$, {\it i.e.}, iff $\beta=2$,
Dyson's model is realized as the Brownian motions
{\it conditioned never to collide with each other}
\cite{Gra99,KT07}.

In order to explain the importance of the notion
of {\it entire functions} for the present problem,
we rewrite the results reported in
our previous paper \cite{KT08} for Dyson's
model with symmetric initial configurations below.
Then the changes which we have to do for
the systems with asymmetric initial configurations
are shown.
There the origin of the quadratic term $t^2/4$
in (\ref{eqn:finite_Airy}) will be clarified.

\subsection{Processes with Space-time
Determinantal Correlations
and Entire Functions}

In an earlier paper \cite{KT08}, we studied a class of
a non-equilibrium dynamics of 
Dyson's model with $\beta=2$ and
an infinite number of particles.
As an example in the class, we reported a 
relaxation process, denoted here by $(\Xi(t), \P_{\sin})$,
which starts from a configuration
\begin{equation}
\xi_{\Z}(\cdot)=\sum_{a \in \Z} \delta_a(\cdot),
\label{eqn:xiZ}
\end{equation}
in which every point of $\Z$ is occupied by one particle, 
and converges to the stationary state $\mu_{\sin}$.
This process $(\Xi(t), \P_{\sin})$ is 
{\it determinantal}, in the sense that
there is a function $\mbK_{\sin}(s,x ; t, y)$
called the {\it correlation kernel} such that
it is continuous with respect to $(x, y) \in \R^2$
for any fixed $(s,t) \in [0, \infty)^2$,
and that, 
for any integer $M \geq 1$,
any sequence $(N_m)_{m=1}^{M}$ of positive integers,
and any time sequence $0 < t_1 < \cdots < t_M < \infty$,
the $(N_1, \dots, N_M)$-{\it multitime correlation function}
$\rho_{\sin}(t_1, \x^{(1)}_{N_1}; \dots; t_M, \x^{(M)}_{N_M}),
\x^{(m)}_{N_m}= (x^{(m)}_1, \dots, x^{(m)}_{N_m} ) \in \R^{N_m},
1 \leq m \leq M$, is expressed by a determinant of 
a $\sum_{m=1}^{M} N_m \times \sum_{m=1}^{M}N_m$
asymmetric real matrix;
\begin{equation}
\rho_{\sin} \Big( t_1, \x^{(1)}_{N_1}; \dots ; 
t_M, \x^{(M)}_{N_M} \Big)=
\det_{
\substack{1 \leq j \leq N_{m}, 
1 \leq k \leq N_{n} \\ 1 \leq m, n \leq M}
}
\Big[
\mbK_{\sin} (t_m, x_{j}^{(m)}; t_n, x_{k}^{(n)} )
\Big].
\label{eqn:rho0}
\end{equation}
The finite dimensional distributions
of the process $(\Xi(t), \P_{\sin})$ are 
determined by $\mbK_{\sin}$ through (\ref{eqn:rho0}).
It is expected that the correlation kernel $\mbK_{\sin}$ is 
described by using the sine function as is the correlation
kernel $K_{\sin}$ of the stationary distribution
$\mu_{\sin}$ given by (\ref{eqn:Ksin}).
It is indeed true. Set
\begin{equation}
f(z)= \sin (\pi z), \quad z \in \C,
\label{eqn:sine1}
\end{equation}
and 
\begin{equation}
p_{\sin}(t, x)=
\frac{e^{-x^2/2t}}{\sqrt{2 \pi |t|}}, \quad
t \in \R \setminus \{0\}, \quad x \in \C.
\label{eqn:psin}
\end{equation}
When $t >0$, $p_{\sin}(t,y-x)$ is the heat kernel:
the solution of the heat equation
$\partial u(t,x)/\partial t
= (1/2) \partial^2 u(t,x) /\partial x^2$
with $\lim_{t \to 0} u(t,x)dx = \delta_y(dx)$,
and is expressed using (\ref{eqn:sine1}) as
$$
p_{\sin}(t, y-x) =
\frac{1}{2} \int_{\R} du \,
e^{-\pi^2 u^2 t/2} 
\Big\{ f(ux) f(uy)
+ f(ux+1/2) f(uy+1/2) \Big\}.
$$
For $0 < s < t$, by setting (\ref{eqn:psin}),
the Chapman-Kolmogorov equation
\begin{equation}
\int_{\R}dy \, p_{\sin}(t-s, z-y) p_{\sin}(s, y-x)
=p_{\sin}(t, z-x)
\label{eqn:CK1}
\end{equation}
can be extended to
\begin{equation}
\int_{\R}dy \, p_{\sin}(-t, z-y) 
p_{\sin}(t-s, y-x) 
= p_{\sin}(-s, z-x).
\label{eqn:CK2}
\end{equation}
Then $\mbK_{\sin}(s,x;t,y)$ is given by 
\begin{eqnarray}
\mbK_{f}(s,x; t,y)
&=& \sum_{a \in f^{-1}(0)} \int_{\sqrt{-1} \, \R} 
\frac{dz}{\sqrt{-1}} \,
p_f(0,a; s,x) 
\frac{1}{z-a} \frac{f(z)}{f'(a)}
p_f(t, y; 0, z)
\nonumber\\
&& 
- {\bf 1}(s>t) p_f(t, y; s, x),
\qquad s,t \geq 0, \quad x, y \in \R
\label{eqn:mbKf}
\end{eqnarray}
with setting (\ref{eqn:sine1}) and
$p_f(s,x;t,y)=p_{\sin}(t-s, y-x)$ 
with (\ref{eqn:psin}), 
where $f^{-1}(0)$ denotes the 
{\it zero set} of the function $f$;
$f^{-1}(0)=\{z: f(z)=0\}$, 
$f'(a)=df(z)/dz|_{z=a}$, 
and ${\bf 1}(\omega)$ is the indicator of
a condition $\omega$;
${\bf 1}(\omega)=1$ if $\omega$ is satisfied,
and ${\bf 1}(\omega)=0$ otherwise.
In this paper 
$\int_{\sqrt{-1} \, \R} dz \, \cdot$
means the integral on the imaginary axis in $\C$
from $-\sqrt{-1} \infty$ to $\sqrt{-1} \infty$.
The well-definedness of the correlation kernel $\mbK_{\sin}$
and thus of the process $(\Xi(t), \P_{\sin})$
is guaranteed \cite{KT08} by the fact that 
the sine function (\ref{eqn:sine1}) is 
an {\it entire function} ({\it i.e.}, analytic
in the whole complex plane $\C$),
and the {\it order of growth} $\rho_f$, which is
generally defined for an entire function $f$ by
$$
\rho_f= \limsup_{r \to \infty}
\frac{\log \log M_f(r)}{\log r}
\quad \mbox{for} \quad
M_f(r)=\max_{|z|=r} |f(z)|,
$$
is one. 
(The {\it type} defined by
$
\sigma_{f}=\limsup_{r \to \infty} 
\log M_f(r) / r^{\rho_f}
$
is equal to $\pi$ for (\ref{eqn:sine1}).
That is, the sine function (\ref{eqn:sine1})
is an {\it entire function of exponential type}
$\pi$ \cite{Lev96} ;
$
M_f(r) \sim e^{\pi r}
$
as $r \to \infty$.)
We can show that
$$
\mbK_{\sin}(t, x; t, y)\mbK_{\sin}(t, y; t, x) dx dy 
\to
\xi_{\Z}(dx) {\bf 1}(x=y) \qquad
\mbox{as} \quad t \to 0,
$$
since $f^{-1}(0)=\sin^{-1}(0)/\pi=\Z$.
It implies that the initial configuration (\ref{eqn:xiZ})
of the relaxation process
$(\Xi(t), \P_{\sin})$ shall be regarded as the point-mass
distribution on the zero set of the sine function (\ref{eqn:sine1}).
Moreover, we showed in \cite{KT08}, by noting
$$
\frac{1}{z-a} \frac{f(z)}{f'(a)}
=K_{\sin}(z-a),
$$
if $a \in f^{-1}(0)=\Z$ and
$z \not=a$ for (\ref{eqn:sine1}),
that 
\begin{equation}
\mbK_{\sin}(s+\theta, x; t + \theta, y)
\to
\bK_{\sin}(t-s, y-x)
\quad \mbox{as} \quad \theta \to \infty
\label{eqn:convKsine}
\end{equation}
with the so-called {\it extended sine kernel},
\begin{equation}
\bK_{\sin}(t, x) 
= \left\{ \begin{array}{ll} 
\displaystyle{
\int_{0}^{1} du \, e^{\pi^2 u^2 t/2} 
\cos ( \pi u x ) }
& \mbox{if $t>0 $} \cr
& \cr
K_{\sin}(x)
& \mbox{if $t=0$} \cr
& \cr
\displaystyle{
- \int_{1}^{\infty} du \, 
e^{\pi^2 u^2 t/2} \cos (\pi u x)}
& \mbox{if $t<0$},
\end{array} \right.
\label{eqn:sine-kernel}
\end{equation}
$x \in \R$.
The equilibrium dynamics in $\mu_{\sin}$,
first studied by Spohn \cite{Spo87},
has been shown to be 
determinantal with the correlation kernel 
(\ref{eqn:sine-kernel}) by
Nagao and Forrester \cite{NF98}.
This process is realized in
the long-term limit of the relaxation process
$(\Xi(t), \P_{\sin})$.
See also the Dirichlet form approach by Osada to
the reversible process with respect to 
$\mu_{\sin}$ \cite{Osa96,Osa08}.

Now we set
\begin{equation}
  f(z)=\Ai(z), \quad z \in \C.
\label{eqn:f=Ai}
\end{equation}
The Airy function $\Ai(z)$, (\ref{eqn:Airy1}), 
is another entire function, 
whose order of growth is $\rho_f=3/2$ with
type $\sigma_f=2/3$;
$
\max_{|z|=r} |f(z)| \sim \exp [ (2/3) r^{3/2} ]
$
as $r \to \infty$.
For $t>0$, we consider
\begin{equation}
p_{\Ai}(t, y|x)=\int_{\R} du \,
e^{ut/2} f(u+x) f(u+y), \quad
x, y \in \R,
\label{eqn:pAi0}
\end{equation}
which is the solution 
of the differential equation; 
$$
\frac{\partial}{\partial t} u(t,x)
= \frac{1}{2} \left(
\frac{\partial^2}{\partial x^2} - x \right) u(t,x)
\quad \mbox{with} \quad
\lim_{t \to 0} u(t, x)dx= \delta_y (dx).
$$
The integral
$\displaystyle{\int_{\R} dz \, p_{\Ai}(t, z|x)}$
is given by 
\begin{equation}
g(t,x) =
\exp \left( - \frac{tx}{2}+\frac{t^3}{24} \right).
\label{eqn:g1}
\end{equation}
We find, for $s < t < 0$,
$g(s, x) p_{\Ai}(t-s,y|x)/g(t, y)$
is equal to the transition probability density 
of 
\begin{equation}
B(t)+\frac{t^2}{4}
\label{eqn:tsquare}
\end{equation}
from $x$ at time $s$ to $y$ at time $t$,
$x, y \in \R$, 
where $B(t), t \in [0, \infty)$ is the one-dimensional
standard Brownian motion.
(See also \cite{FS05} and references therein.)
Then for $s, t \in \R, s \not= t, x, y \in \C$, 
we set
\begin{eqnarray}
&& q(s, t, y-x) =
p_{\sin}\left(t-s, \left(y-\frac{t^2}{4} \right)-
\left( x-\frac{s^2}{4} \right) \right)
\nonumber\\
\label{eqn:q}
&& = \frac{1}{\sqrt{2 \pi |t-s|}}
\exp \left[ 
-\frac{(y-x)^2}{2(t-s)}+\frac{(t+s)(y-x)}{4}
-\frac{(t-s)(t+s)^2}{32} \right],
\end{eqnarray}
and as an extension of (\ref{eqn:pAi0}) we define
\begin{eqnarray}
&& p_{\Ai}(t-s, y|x)
= \frac{g(t,y)}{g(s,x)} q(s, t, y-x) 
\nonumber\\
\label{eqn:pAi}
&& = \frac{1}{\sqrt{2\pi|t-s|}}
\exp \left[ -\frac{(y-x)^2}{2(t-s)}
-\frac{(t-s)(y+x)}{4} + \frac{(t-s)^3}{96} \right].
\end{eqnarray}
Corresponding to (\ref{eqn:CK1}) and (\ref{eqn:CK2}),
we have the two sets of equalities
\begin{eqnarray}
\label{eqn:CKq1}
&& \int_{\R}dy \, q(s,t, z-y) q(0,s,y-x)=q(0,t,z-x) \\
\label{eqn:CKq2}
&& \int_{\R} dy \, q(t,0,z-y) q(s,t, y-x)=q(s,0,z-x)
\end{eqnarray}
and
\begin{eqnarray}
\label{eqn:CKAi1}
&& \int_{\R} dy \, p_{\Ai}(t-s,z|y) p_{\Ai}(s,y|x)
=p_{\Ai}(t, z|x) \\
\label{eqn:CKAi2}
&& \int_{\R} dy \,
p_{\Ai}(-t,z|y) p_{\Ai}(t-s,y|x)
= p_{\Ai}(-s,z|x)
\end{eqnarray}
for $0 < s < t$.

Let $\mbK_{\Ai}$ be the function
given by (\ref{eqn:mbKf}) with setting
(\ref{eqn:f=Ai}) and
$p_f(s,x;t,y)=p_{\Ai}(t-s,y|x)$
with (\ref{eqn:pAi}).
We will prove that
$\mbK_{\Ai}$ is well-defined as a correlation kernel and
it determines finite dimensional distributions of
an infinite particle system
through a similar formula to (\ref{eqn:rho0}).
We denote this system by $(\Xi_{\cA}(t), \P_{\Ai})$. 
The fact that $p_{\Ai}$ used in $\mbK_{\Ai}$ is 
a transform (\ref{eqn:pAi}) of the transition probability
density $q$ of (\ref{eqn:tsquare}) is the origin
of the quadratic term $t^2/4$ in (\ref{eqn:finite_Airy}).
We can show that
$$
\mbK_{\Ai}(t, x; t, y) \mbK_{\Ai}(t, y; t, x) dx dy 
\to
\xi_{\cA}(dx) {\bf 1}(x=y) \quad \mbox{as} \quad
t \to 0.
$$
By using the integral formula 
for (\ref{eqn:f=Ai})
$$
\frac{1}{z-a} \frac{f(z)}{f'(a)}
=\frac{1}{(\Ai'(a))^2}
\int_{0}^{\infty} du \, \Ai(u+z) \Ai(u+a)
$$
for $a \in {\cA}, z \not= a$,
and the fact that
$
\{ \Ai(x+a)/\Ai'(a), 
a \in {\cA} \}
$
forms a complete orthonormal basis 
for the space $L^2(0, \infty)$ of
square integrable functions 
on the interval $(0, \infty)$ \cite{Tit62}, 
we will prove that
\begin{equation}
\mbK_{\Ai}(s+\theta, x; t + \theta, y)
\to
{\bK}_{\Ai}(t-s, y|x)
\quad \mbox{as} \quad \theta \to \infty,
\label{eqn:convKAi}
\end{equation}
where ${\bK}_{\Ai}$ is the
so-called {\it extended Airy kernel},
\begin{equation}
{\bK}_{\Ai}(t, y|x)
= \left\{ 
\begin{array}{cc}
\displaystyle{
\int_{0}^{\infty} d u \,
e^{-ut/2} \Ai(u+x) \Ai(u+y)}
& \mbox{if $t \geq 0$} \cr
& \cr
\displaystyle{
- \int_{-\infty}^{0} d u \,
e^{-ut/2} \Ai(u+x) \Ai(u+y) }
& \mbox{if $t < 0$},
\end{array}
\right.
\label{eqn:AiryK}
\end{equation}
$x, y \in \R$.
We denote by $(\Xi_{\cA}(t), {\bP}_{\Ai})$ 
the infinite particle system, which is determinantal
with the correlation kernel ${\bK}_{\Ai}$
\cite{FNH99,NKT03,KNT04,KT07}.
The Airy kernel (\ref{eqn:AiryKer}) of $\mu_{\Ai}$ is
given by $K_{\Ai}(y|x)={\bf K}_{\Ai}(0, y|x)$ and thus 
$(\Xi_{\cA}(t), {\bP}_{\Ai})$ is a reversible process with respect to
$\mu_{\Ai}$. The process $(\Xi_{\cA}(t), \P_{\Ai})$,
which is determinantal with the correlation kernel
$\mbK_{\Ai}$,
is a non-equilibrium infinite particle system
exhibiting the relaxation phenomenon from the initial configuration
$\xi_{\cA}$ to the stationary state $\mu_{\Ai}$.

Then consider the finite-particle system (\ref{eqn:finite_Airy}) again.
Let $\P_{\cA}^{\xi^N}$ be the distribution of 
the process $\Xi_{\cA}(t)=\sum_{j=1}^{N} \delta_{Y_j(t)}$
starting from a configuration $\xi^N$.
We denote by $\mM$ the space of nonnegative 
integer-valued Radon measures on $\R$,
which is a Polish space with the {\it vague topology}:
we say $\xi_n$ converges to $\xi$ vaguely, if 
$\displaystyle{\lim_{n \to \infty} \int_{\R} \varphi(x) \xi_n(dx)
=\int_{\R} \varphi(x) \xi(dx)}$ 
for any $\varphi \in {\rm C}_0(\R)$,
where ${\rm C}_0(\R)$ is the set of all 
continuous real-valued functions with
compact supports.
Any element $\xi$ of $\mM$ can be represented as
$\xi(\cdot) = \sum_{j\in \Lambda}\delta_{x_j}(\cdot)$
with an index set $\Lambda$
and a sequence of points in $\R$, $\x =(x_j)_{j \in \Lambda}$ 
satisfying $\xi(I)=\sharp\{x_j: x_j \in I \} < \infty$ 
for any compact subset $I \subset \R$.
For $A \subset \R$, we write the {\it restriction} of
$\xi$ on $A$ as
$(\xi\cap A) (\cdot)= \displaystyle{\sum_{j \in \Lambda : x_j \in  A}}
\delta_{x_j}(\cdot)$.
We put
$
\mM_0= \Big\{ \xi\in\mM : 
\xi(\{x\})\le 1 \mbox { for any }  x\in\R\Big\}.
$
We will prove that the finite particle process
$(\Xi_{\cA}(t), \P_{\cA}^{\xi^N})$ is determinantal
for any initial configuration $\xi^N \in \mM_0$
and give the correlation kernel $\mbK_{\cA}^{\xi^N}$
(Proposition \ref{Proposition:Finite}).
For $\xi \in \mM$ with an infinite number 
of particles $\xi(\R)=\infty$, 
when $\mbK_{\cA}^{\xi \cap [-L, L]}$
converges to a continuous function as $L \to \infty$,
the limit is written as $\mbK_{\cA}^{\xi}$. 
If $\P^{\xi \cap [-L, L]}_{\cA}$ 
converges to a probability measure $\P^{\xi}_{\cA}$
on $\mM^{[0, \infty)}$,
which is determinantal with the correlation kernel
$\mbK_{\cA}^{\xi}$, weakly in the sense of finite dimensional
distributions as $L \to \infty$ in the vague topology, 
we say that the process $(\Xi_{\cA}(t), \P_{\cA}^{\xi})$
is {\it well defined with the correlation kernel}
$\mbK_{\cA}^{\xi}$.
(The regularity of the sample paths 
of $\Xi_{\cA}(t)$ will be discussed 
in the forthcoming paper \cite{KT09}.)
We will give sufficient conditions for initial 
configurations $\xi \in \mM_0$ so that
the process $(\Xi_{\cA}(t), \P_{\cA}^{\xi})$
is well defined (Theorem \ref{Theorem:Infinite}). 
We denote by $\mX^{\cA}$ the set of configurations $\xi$
satisfying the conditions and put 
$\mX^{\cA}_0=\mX^{\cA} \cap \mM_0$.
It is clear that the configuration $\xi_{\cA} \in \mX^{\cA}_0$.
Then, if we consider the finite particle systems
$\Xi_{\cA}(t)=\sum_{j=1}^{N} \delta_{Y_j(t)}, N \geq 2$, 
with (\ref{eqn:finite_Airy}) 
starting from the $N$-particle approximation of $\xi_{\cA}$, 
(\ref{eqn:xiAN}),
we can prove
$(\Xi_{\cA}(t), \P_{\cA}^{\xi_{\cA}^N})
\to (\Xi_{\cA}(t), \P_{\Ai})$ as $N \to \infty$
in the sense of finite dimensional distributions
(Theorem \ref{Theorem:from_roots} (i)).
That is, 
$(\Xi_{\cA}(t), \P_{\Ai})=(\Xi_{\cA}(t), \P_{\cA}^{\xi_{\cA}})$ 
with (\ref{eqn:xiA}).
Moreover, we will show (\ref{eqn:convKAi})
and prove the relaxation phenomenon
$(\Xi_{\cA}(t+\theta), \P_{\Ai})
\to (\Xi_{\cA}(t), \bP_{\Ai})$ 
(Theorem \ref{Theorem:from_roots} (ii)).

The paper is organized as follows.
In Sect. 2 preliminaries and main results are given.
Some remarks on extensions of the present results
are also given there.
In Sect. 3 the properties of the Airy function
used in this paper are summarized.
Section 4 is devoted to proofs of results.

\SSC{Preliminaries and Main Results}

For $\xi(\cdot)=\sum_{j\in \Lambda}
\delta_{x_j}(\cdot) \in\mM$,
we introduce the following operations;
\begin{description}


\item[(shift)] for $u \in \R$, 
$\tau_u \xi(\cdot) =\displaystyle{\sum_{j \in \Lambda}} 
\delta_{x_j+u}(\cdot)$,


\item[(square)]
$\displaystyle{
\xi^{\langle 2 \rangle}(\cdot)
=\sum_{j \in \Lambda} \delta_{x_j^2} (\cdot)}$.
\end{description}
We use the convention such that
$$
\prod_{x\in\xi}f(x) =\exp
\left\{\int_\R \xi(dx) \log f(x) \right\}
=\prod_{x \in \supp \xi}f(x)^{\xi(\{x\})}
$$
for $\xi\in \mM$ and a function $f$ on $\R$,
where $\supp \xi = \{x \in \R : \xi(\{x\}) > 0\}$.
For a multivariate symmetric function $g$ we write 
$g((x)_{x \in \xi})$ for $g((x_j)_{j \in \Lambda})$.

\subsection{Determinantal processes}

As an $\mM$-valued process $(\Xi(t), \P^{\xi})$,
we consider the system such that,
for any integer $M \geq 1$,
$f_m \in {\rm C}_{0}(\R), \theta_m \in \R,
1 \leq m \leq M$,
$0 < t_1 < \cdots < t_M < \infty$,
the expectation of
$\displaystyle{\exp \Big\{ \sum_{m=1}^{M} \theta_m 
\int_{\R} f_m(x) \Xi(t_m, dx) \Big\}}$
can be expanded by
$$
\chi_{m}(x)=e^{\theta_{m} f_{m}(x)}-1, \quad
1 \leq m \leq M, 
$$
as
\begin{eqnarray}
{\cG}^{\xi}[\chi] &\equiv& 
\E^{\xi} \left[
\exp \left\{ \sum_{m=1}^{M} \theta_m 
\int_{\R} f_m(x) \Xi(t_m, dx) \right\} \right]
\nonumber\\
&=& 
\sum_{N_{1} \geq 0} \cdots
\sum_{N_{M} \geq 0}
\prod_{m=1}^{M}\frac{1}{N_m !}
\int_{\R^{N_{1}}} \prod_{j=1}^{N_1} d x_{j}^{(1)}
 \cdots
\int_{\R^{N_{M}}} 
\prod_{j=1}^{N_{M}} d x_{j}^{(M)} \nonumber\\
&& \qquad \times \prod_{m=1}^{M} \prod_{j=1}^{N_{m}} 
\chi_{m} \Big(x_{j}^{(m)} \Big)
\rho\Big( t_{1}, \x^{(1)}; \dots ; t_{M}, \x^{(M)} \Big).
\label{eqn:Gxi}
\end{eqnarray}
Here $\rho$'s are locally integrable functions,
which are symmetric in the sense that
$\rho(\dots; t_m, \sigma(\x^{(m)}); \dots)
=\rho(\dots; t_m, \x^{(m)}; \dots)$
with $\sigma(\x^{(m)})
\equiv (x^{(m)}_{\sigma(1)}, \dots, x^{(m)}_{\sigma(N_m)})$
for any permutation $\sigma \in {\cS}_{N_m},
1 \leq \forall m \leq M$.
In such a system $\rho( t_{1}, \x^{(1)};
\dots ; t_{M}, \x^{(M)})$ is called
the $(N_1, \dots, N_{M})$-{\it multitime correlation function}
and ${\cG}^{\xi}[\chi]$ the {\it generating function 
of multitime correlation functions}.
There are no multiple points with probability one for $t > 0$.
Then we assume that there is a function $\mbK(s,x;t,y)$,
which is continuous with respect to
$(x,y) \in \R^2$ for any fixed $(s,t) \in [0, \infty)^2$,
such that
$$
\rho \Big(t_1,\x^{(1)}; \dots;t_M,\x^{(M)} \Big) 
=\det_{
\substack{1 \leq j \leq N_{m}, 1 \leq k \leq N_{n} 
\\ 1 \leq m, n \leq M}
}
\Bigg[
\mbK(t_m, x_{j}^{(m)}; t_n, x_{k}^{(n)} )
\Bigg]
$$
for any integer $M \geq 1$, 
any sequence $(N_m)_{m=1}^{M}$ of positive integers, and
any time sequence $0 < t_1 < \cdots < t_M < \infty$.
Let $\T=\{t_1, \dots, t_M\}$.
We note that
$\Xi^{\T}=\sum_{t \in \T} \delta_t \otimes \Xi(t)$
is a determinantal (Fermion) point process on $\T \times \R$
with an operator ${\cK}$ given by
$$
{\cK}f(s,x)=\sum_{t \in \T} \int_{\R} dy \,
\mbK(s,x;t, y) f(t,y),
\quad f(t, \cdot) \in {\rm C}_0(\R),\,  t \in \T.
$$
When ${\cK}$ is symmetric, Soshnikov \cite{Sos00}
and Shirai and Takahashi \cite{ST03} gave sufficient
conditions for $\mbK$ to be a correlation kernel
of a determinantal point process.
Though such conditions are not known for asymmetric cases,
a variety of processes, which are determinantal
with asymmetric correlation kernels,
have been studied. See, for example, \cite{TW03,KT07}.
If there exists a function $\mbK$,
which has the above properties and determines
the finite dimensional distributions of the
process $(\Xi(t), \P^{\xi})$, we say the process 
$(\Xi(t), \P^{\xi})$ is {\it determinantal with the 
correlation kernel} $\mbK$ \cite{KT08}.

For $N \in \N$, the determinant of an $N \times N$ matrix
${\rm M}=(m_{jk})_{1 \leq j, k \leq N}$
is defined by
$\displaystyle{\sum_{\sigma \in {\cS}_N}
{\rm sgn}(\sigma) \prod_{j=1}^{N} m_{j \sigma(j)}}$,
where ${\rm sgn}(\sigma)$ denotes the sign of 
permutation $\sigma$. Any permutation $\sigma$
consists of exclusive cycles.
If we write each cyclic permutation as
$$
{\sf c}=\left( \begin{array}{cccc}
a & b & \cdots & \omega \cr
b & c & \cdots & a \end{array} \right)
$$
and the number of cyclic permutations in a given $\sigma$
as $\ell(\sigma)$, then the determinant of ${\rm M}$ is
expressed as
$$
\det {\rm M}=\sum_{\sigma \in {\cS}_N}
{\rm sgn}(\sigma) 
\prod_{{\sf c}_j: 1 \leq j \leq \ell(\sigma)}
\Big( m_{ab} m_{bc} \dots m_{\omega a} \Big).
$$
It implies that, with given $a_1, a_2, \dots, a_N$,
even if each element $m_{jk}$ of the matrix ${\rm M}$
is replaced by $m_{jk} \times (a_j/a_k)$,
the value of determinant is not changed.
The above observation will lead to the following lemma.

\begin{lem}\label{thm:gauge}
Let $(\Xi(t), \P)$ and $(\widetilde{\Xi}(t), \widetilde{\P})$
be the processes, which are determinantal with
correlation kernels $\mbK$ and $\widetilde{\mbK}$,
respectively.
If there is a function $G(s,x)$, which is continuous 
with respect to $x \in \R$ for any fixed $s \in [0, \infty)$,
such that
\begin{equation}
\mbK(s,x;t,y)=\frac{G(s,x)}{G(t,y)}
\widetilde{\mbK}(s,x;t,y),
\quad s, t \in [0, \infty), \quad x, y \in \R,
\label{eqn:gauge1}
\end{equation}
then
\begin{equation}
(\Xi(t), \P)=(\widetilde{\Xi}(t), \widetilde{\P})
\label{eqn:gauge2}
\end{equation}
in the sense of finite dimensional distributions.
\end{lem}
\vskip 0.3cm

In literatures, (\ref{eqn:gauge1})
is called the {\it gauge transformation} and
(\ref{eqn:gauge2}) is said to be the {\it gauge invariance}
of the determinantal processes.

\subsection{The Weierstrass canonical product
and entire functions}

For $\xi^N \in \mM_0, \xi^N(\R)=N < \infty$,
with $p \in \N_0 \equiv \N \cup \{0\}$
we consider the product
$$
\Pi_{p}(\xi^N, z)= 
\prod_{x \in \xi^{N} \cap \{0\}^{\rm c}}
G \left( \frac{z}{x}, p \right),
\quad z \in \C,
$$
where 
\begin{equation}
G(u,p) = \left\{
   \begin{array}{ll}
\displaystyle{1-u} 
& \mbox{if} \quad p=0  \\
& \\
\displaystyle{(1-u)\exp
\left [ u+\frac{u^2}{2}+\cdots +\frac{u^p}{p} \right]}
& \mbox{if} \quad p\in\N.
   \end{array} \right. 
\label{Weierstrass}
\end{equation}
The functions $G(u,p)$ are called
the {\it Weierstrass primary factors}.
With $\alpha > 0$ we put
$$
M_\alpha(\xi^N)
=\left( \int_{\{0\}^{\rm c}} 
\frac{1}{|x|^\alpha}\xi^N(dx)
\right)^{1/\alpha}.
$$
For $\xi \in \mM_0$ with $\xi(\R)=\infty$,
we write 
$M_\alpha(\xi,L)$ for $M_\alpha(\xi\cap [-L,L]), L > 0$, 
and put
$\displaystyle{
M_\alpha(\xi)= \lim_{L\to\infty}M_\alpha(\xi,L),
}$
if the limit finitely exists.
If $M_{p+1}(\xi) < \infty$ for some $p \in \N_0$,
the limit
\begin{equation}
\Pi_{p}(\xi, z)= \lim_{L \to \infty} 
\Pi_{p}(\xi \cap [-L, L], z)
= \prod_{x \in \xi \cap \{0\}^{\rm c}}
G \left( \frac{z}{x}, p \right),
\quad z \in \C
\label{eqn:Pip}
\end{equation}
finitely exists.
This infinite product is called the
{\it Weierstrass canonical product of genus} $p$ \cite{Lev96}.
The {\it Hadamard theorem} \cite{Lev96} claims that any entire function $f$
of finite order $\rho_f < \infty$ can be represented by
\begin{equation}
f(z)=z^m e^{P_{q}(z)}
\Pi_{p}(\xi_f, z),
\label{eqn:Hadamard}
\end{equation}
where $p$ is a nonnegative integer less than or equal 
to $\rho_f$,
$P_q(z)$ is a polynomial in $z$ of degree $q\le \rho_f$, 
$m$ is the multiplicity of the root at the origin,
and 
$
\xi_f=\sum_{x \in f^{-1}(0) \cap \{0\}^{\rm c}}
\delta_x.
$
We give two examples ;
\begin{eqnarray}
\label{eqn:WH1a}
\sin( \pi z) &=& \pi z \Pi_0(\xi_{\Z}, z), \\
\label{eqn:WH1}
\Ai(z) &=& e^{d_0+d_1 z} \Pi_1(\xi_{\cA}, z)
\end{eqnarray}
with (\ref{eqn:xiZ}), (\ref{eqn:xiA}) and
\begin{eqnarray}
d_0 &=&\log \Ai (0) = - \log \Big( 3^{2/3} \Gamma(2/3) \Big),
\nonumber\\
d_1 &=& \frac{\Ai'(0)}{\Ai(0)}
=-\frac{3^{1/3} \Gamma(2/3)}{\Gamma(1/3)}
=- \frac{3^{5/6} (\Gamma(2/3))^2}{2 \pi}.
\label{eqn:d0d1}
\end{eqnarray}

For $\xi^N\in\mM_0$ with $\xi^N(\R)=N$ we put
\begin{equation}
\Phi_p(\xi^N, a, z) \equiv \Pi_p(\tau_{-a}\xi^N, z-a)
=\prod_{x\in\xi^N \cap \{a\}^{\rm c}}
G\left(\frac{z-a}{x-a},p \right),
\quad a, z \in \C.
\label{entire2}
\end{equation}
With (\ref{eqn:xiAN})
we set
\begin{eqnarray}
\Phi_{\cA}(\xi^N,z) 
&\equiv& e^{d_1 z} \exp \Bigg[
\int_{\R} \frac{z}{x} \xi_{\cA}^N(dx) \Bigg]
\Pi_0(\xi^N, z)
\nonumber\\
\label{eqn:entire_A0}
&=& e^{d_1 z}
\exp \Bigg[
\int_{\{0\}^{\rm c}} \frac{z}{x} (\xi_{\cA}^N-\xi^N)(dx)
\Bigg]
\Pi_1(\xi^N, z), \quad z \in \C, \\
\Phi_{\cA}(\xi^N,a, z) &\equiv& \Phi_{\cA}(\tau_{-a}\xi^N,z-a)
\nonumber\\
\label{eqn:entire_A} 
&=& e^{d_1(z-a)} 
\exp \Bigg[
\int_{\R} \frac{z-a}{x} \xi_{\cA}^N(dx) \Bigg]
\Phi_0(\xi^N, a, z),
\quad a, z \in \C. \nonumber\\
\end{eqnarray}

\begin{lem}\label{thm:2_1}
Let $\xi^N\in\mM_0$ with $\xi^N(\R)=N < \infty$ and
$\xi^N(\{0\})=0$.
Then for $a \in \supp \xi^{N}, z \not= a$,
\begin{equation}
\Phi_{\cA}(\xi^N,a,z)
= \frac{1}{z-a}\frac{\Phi_{\cA}(\xi^N,z)}{\Phi_{\cA}'(\xi^N,a)},
\label{eqn:relation2_A}
\end{equation}
where $\Phi_{\cA}'(\cdot,a)=\partial \Phi_{\cA}(\cdot,z)/
\partial z \Big|_{z=a}$.
\end{lem}

\vskip 3mm
\noindent {\it Proof.}
Since
$1- (z-a)/(x-a) = (x-z)/(x-a)
= (1-z/x)/(1-a/x)$,
\begin{equation}
\Phi_{\cA}(\xi^N,a,z)=
\frac{\displaystyle{e^{d_1 z} \exp \left[
\int_{\R} \frac{z}{x} \xi^{N}_{\cA}(dx) \right]
\Pi_0(\xi^{N}-\delta_a,z)}}
{\displaystyle{e^{d_1 a} \exp \left[
\int_{\R} \frac{a}{x} \xi^{N}_{\cA}(dx) \right]
\Pi_0(\xi^{N}-\delta_a,a)}},
\label{eqn:eqQ}
\end{equation}
where the numerator is equal to
$\Phi_{\cA}(\xi^N,z)/(1-z/a)$.
From (\ref{eqn:entire_A0}), we have
\begin{eqnarray}
&& \frac{\partial}{\partial z} 
\Phi_{\cA}(\xi^N, z) = d_1 \Phi_{\cA}(\xi^N, z) \nonumber\\
&& \qquad + e^{d_1 z} 
\int_{\R} \frac{1}{x} \xi_{\cA}^{N}(dx) 
\exp \left[ \int_{\R} \frac{z}{y} \xi_{\cA}^{N}(dy) \right]
\Pi_{0}(\xi^N, z)
\nonumber\\
&& \qquad + e^{d_1 z}  
\exp \left[ \int_{\R} \frac{z}{y} \xi_{\cA}^{N}(dy) \right]
\int_{\R} \left(-\frac{1}{x}\right) 
\Pi_{0}(\xi^N-\delta_x, z) \xi^{N}(dx).
\nonumber
\end{eqnarray}
Since $a \in \supp \xi^N$ is assumed,
$$
\Phi'_{\cA}(\xi^N, a)=
e^{d_1 a} \exp \left[ \int_{\R} \frac{a}{x} \xi^{N}(dx) \right]
\left(-\frac{1}{a}\right)
\Pi_{0}(\xi^N-\delta_a, a).
$$
It implies that the denominator of (\ref{eqn:eqQ})
is equal to $-a \Phi'_{\cA}(\xi^N, a)$.
Then (\ref{eqn:relation2_A}) is obtained.
\qed
\vskip 0.5cm

For $\xi^N\in\mM_0$ with $\xi^N(\R)=N$ we put
\begin{equation}
M_{\cA}(\xi^N)
=\int_{\{0\}^{\rm c}} \frac{1}{x} (\xi_{\cA}^N-\xi^N)(dx).
\label{eqn:MAN}
\end{equation}
For $\xi\in\mM_0$ with $\xi(\R)=\infty$
we write 
$M_{\cA}(\xi, L)$ for $M_{\cA}(\xi \cap [-L,L]), L > 0$, 
and put
$\displaystyle{
M_{\cA}(\xi)=\lim_{L\to\infty}M_{\cA}(\xi,L),
}$
if the limit finitely exists.
For $\xi \in \mM_0$, 
$p\in\N_0$, $a\in\R$, and $z\in\C$ we define
$\displaystyle{
\Phi_p (\xi,a,z)=\lim_{L\to\infty}
\Phi_p(\xi \cap [a-L, a+L], a, z)
}$
and
$\displaystyle{
\Phi_{\cA} (\xi,a,z)=\lim_{L\to\infty}
\Phi_{\cA}(\xi \cap [a-L, a+L], a, z)
}$,
if the limits finitely exist.
We note that 
$\Phi_p (\xi,a,z)$ finitely exists and
is not identically $0$,
if $M_{p+1}(\tau_{-a}\xi)<\infty$,
and 
$\Phi_{\cA} (\xi,a,z)$ does and 
$\Phi_{\cA}(\xi,a,z) \not \equiv0$,
if $|M_{\cA}(\tau_{-a}\xi)|<\infty$
and $M_2(\tau_{-a} \xi) < \infty$.
For $\xi \in \mM_0, a \in \supp \xi$,
the following equalities will hold, if
all the entries of them
finitely exist ;
\begin{eqnarray}
&& \Phi_0(\xi, a, z) = \Pi_0(\xi, z)
\Phi_0(\xi \cap \{0\}^{\rm c}, a, 0)
\left(\frac{z}{a}\right)^{\xi(\{0 \})}
\frac{a}{a-z},
\nonumber\\
&& \Phi_0(\xi \cap \{0\}^{\rm c},a,0)
=\Pi_0(\xi \cap \{-a\}^{\rm c},-a)
\Phi_0(\xi^{\langle 2 \rangle} \cap \{0\}^{\rm c}, a^2,0)
2^{1-\xi(\{-a\})},
\nonumber
\end{eqnarray}
and then
\begin{eqnarray}
\Phi_1(\xi,a,z) &=& e^{S(\xi,a,z)} \Pi_1(\xi,z)
\Pi_1(\xi\cap \{-a\}^{\rm c},-a) 
\nonumber\\
\label{eqn:finite_zeta}
&& \, \times
\Phi_0(\xi^{\langle 2 \rangle} \cap \{0\}^{\rm c}, a^2,0)
\left(\frac{z}{a}\right)^{\xi(\{0 \})}\frac{a}{a-z},
\end{eqnarray}
where
\begin{equation}
S(\xi,a,z)= \int_{\{a\}^{\rm c}} 
\frac{z-a}{x-a} \xi(dx)
-\int_{\{0\}^{\rm c}} 
\frac{z}{x} \xi(dx)
+ \int_{\{0,-a\}^{\rm c}} 
\frac{a}{x} \xi(dx).
\label{eqn:S1}
\end{equation}

\begin{lem}\label{thm:2_2}
For $a \in \cA, z \not= a$
\begin{equation}
\frac{1}{z-a}
\frac{\Ai(z)}{\Ai'(a)}
= \Phi_{\cA}(\xi_{\cA}, a, z).
\label{eqn:relation1b}
\end{equation}
\end{lem}
\noindent{\it Proof.} 
By (\ref{eqn:WH1}) and the definition (\ref{eqn:entire_A0}),
\begin{equation}
\Ai (z) 
=e^{d_0} \Phi_{\cA}(\xi_{\cA},z), \quad z \in \C.
\label{eqn:relation1}
\end{equation}
As approximations of the Airy function we introduce functions
\begin{equation}
\Ai_N(z) = e^{d_0+d_1 z}
\prod_{\ell=1}^{N} \left( 1- \frac{z}{a_{\ell}} \right)
e^{z/a_{\ell}}, \quad N \in \N,
\label{eqn:WH1app}
\end{equation}
where $0 >a_1 > \cdots > a_N$ are
the first $N$ zeros of $\Ai(z)$.
Since $\xi_{\cA}^N=\sum_{j=1}^{N} \delta_{a_j}$
satisfies the condition of Lemma \ref{thm:2_1},
(\ref{eqn:relation2_A}) with (\ref{eqn:relation1})
and 
$\Ai_N'(a_j)=e^{d_0}
\Phi_{\cA}'(\xi_{\cA}^{N}, a_j)$
gives
$$
\frac{1}{z-a_{j}} 
\frac{\Ai_N(z)}{\Ai_N'(a_j)}
= \Phi_{\cA}(\xi_{\cA}^N, a_j, z),
\quad 1 \leq j \leq N.
$$
Taking $N\to\infty$, we have (\ref{eqn:relation1b}).
\qed

\subsection{Statement of results}

For the solution $\X(t)=(X_1(t),X_2(t),\dots,X_N(t))$ 
of Dyson's model (\ref{eqn:Dyson}) with $\beta=2$
with the initial state $\X(0)=\x$, 
we denote the distribution of the process 
$\Xi(t)=\sum_{j=1}^N \delta_{X_j(t)}$ by $\P^{\xi^N}$
with $\xi^N=\sum_{j=1}^{N} \delta_{x_j}$.
In \cite{KT08} we proved that Dyson's model
(\ref{eqn:Dyson}) with $\beta=2$
starting from any fixed configuration 
$\xi^N \in \mM$ 
is determinantal with the correlation kernel $\mbK^{\xi^N}$
given by
\begin{eqnarray}
&& \mbK^{\xi^N}(s, x; t, y)
= \frac{1}{2 \pi \sqrt{-1}} 
\oint_{\Gamma(\xi^N)} dz
\int_{\sqrt{-1} \, \R} \frac{dw}{\sqrt{-1}}
\, 
\nonumber\\
&& \hskip 3cm \times p_{\sin}(s, x-z)
\frac{1}{w-z}
\prod_{x' \in \xi^{N}}
\left( 1- \frac{w-z}{x'-z} \right)
p_{\sin}(-t, w-y)
\nonumber\\
&&\qquad\qquad\qquad - {\bf 1}(s > t)p_{\sin}(s-t, y-x),
\label{eqn:KN1a}
\end{eqnarray}
where $\Gamma(\xi^N)$ is a closed contour on the
complex plane $\C$ encircling the points in 
$\supp \xi^N$ on $\R$
once in the positive direction,
and $p_{\sin}$ is given by (\ref{eqn:psin}).
If $\xi^N\in \mM_0$,
by performing the Cauchy integrals
(\ref{eqn:KN1a}) is written as
\begin{eqnarray}
\mbK^{\xi^N}(s, x; t, y)&=& \int_{\R} \xi^N(d x') \, 
\int_{\sqrt{-1} \, \R} 
\frac{dy'}{\sqrt{-1}} \, 
p_{\sin}(s, x-x') 
\Phi_0(\xi^{N}, x', y') p_{\sin}(-t, y'-y)
\nonumber\\
&& - {\bf 1}(s > t)p_{\sin}(s-t, x-y).
\label{eqn:KN1}
\end{eqnarray}
Then the following is obtained for 
the process $(\Xi_{\cA}(t), \P_{\cA}^{\xi^N})$
with $\Xi_{\cA}(t)=\sum_{j=1}^{N} \delta_{Y_j(t)}$,
where $\Y(t)=(Y_1(t), \dots, Y_N(t))$ is given by
(\ref{eqn:finite_Airy}).

\begin{prop}
\label{Proposition:Finite}
The process  $(\Xi_{\cA}(t), \P_{\cA}^{\xi^N})$,
starting from any fixed configuration $\xi^N\in \mM_0$ 
with $\xi^N(\R) = N < \infty$,
is determinantal with the correlation kernel $\mbK_{\cA}^{\xi^N}$
given by 
\begin{eqnarray}
\mbK^{\xi^{N}}_{\cA}(s,x ;t,y)
&=&\int_{\R}\xi^N(dx') 
\int_{\sqrt{-1} \, \R} 
\frac{dy'}{\sqrt{-1}} \, 
q(0,s, x-x')
\Phi_{\cA}(\xi^N, x', y')
q(t,0, y'-y)
\nonumber\\
&& - {\bf 1}(s>t)q(t, s, x-y),
\label{eqn:K3}
\end{eqnarray}
where $q$ is given by (\ref{eqn:q}).
\end{prop}
\vskip 0.3cm

We introduce the following conditions:

\vskip 3mm

\noindent ({\bf C.1})
there exists $C_0 > 0$ such that
$|M_{\cA}(\xi)|  < C_0$,

\vskip 3mm

\noindent ({\bf C.2}) (i) 
there exist $\alpha\in (3/2,2)$ and $C_1>0$ such that
$
M_\alpha(\xi) \le C_1,
$ \\
\noindent (ii) 
there exist $\beta >0$ and $C_2 >0$ such that
$$
M_1(\tau_{-a^2} \xi^{\langle 2 \rangle}) \leq C_2
(|a| \vee 1)^{-\beta}
\quad \mbox{for all} \quad a \in \supp \xi.
$$
\vskip 3mm

\noindent We denote by $\mX^{\cA}$ the set of configurations $\xi$
satisfying the conditions ({\bf C.1}) and ({\bf C.2}),
and put $\mX^{\cA}_0 = \mX^{\cA} \cap \mM_0$.

\begin{thm}
\label{Theorem:Infinite}
If $\xi\in \mX^{\cA}_0$, 
the process $(\Xi_{\cA}(t), \P_{\cA}^{\xi})$
is well defined with the correlation kernel
\begin{eqnarray}
&& \mbK^{\xi}_{\cA}(s,x;t,y)
=\int_{\R}\xi(dx') 
\int_{\sqrt{-1} \, \R} 
\frac{dy'}{\sqrt{-1}} \, 
q(0,s, x-x') 
\Phi_{\cA}(\xi, x', y')
q(t, 0, y'-y)
\nonumber\\
&& \qquad\qquad\qquad - {\bf 1}(s>t)q(t, s, x-y).
\label{eqn:Kinfinite}
\end{eqnarray}
\end{thm}
In the proof of this theorem, a useful estimate
of $\Phi_{\cA}$ in (\ref{eqn:Kinfinite}) is 
obtained (Lemma \ref{thm:4_3} {\rm (ii)}).
By virtue of it, we can see
\begin{equation}
\mbK_{\cA}^{\xi}(t, x; t, y)
\mbK_{\cA}^{\xi}(t, y; t, x)
dxdy \to \xi(dx){\bf 1}(x=y)
\quad \mbox{as} \quad t \to 0
\label{eqn:strongC}
\end{equation}
in the vague topology.
Then Theorem \ref{Theorem:Infinite}
gives an infinite particle system
starting form the configuration $\xi$.

The main result of the present paper is the following.

\begin{thm}
\label{Theorem:from_roots}
{\rm (i)} Let $\xi_{\cA}^N(\cdot)$ be the
configuration (\ref{eqn:xiAN}). Then
$$
(\Xi_{\cA}(t), \P_{\cA}^{\xi_{\cA}^N})
\to
(\Xi_{\cA}(t), \P_{\Ai})
\quad \mbox{as} \quad N \to \infty
$$
in the sense of finite dimensional distributions.
Here the process $(\Xi_{\cA}(t), \P_{\Ai})$ is determinantal with
the correlation kernel (\ref{eqn:mbKf}) with setting
(\ref{eqn:f=Ai}) and 
$p_f(s, x; t, y)=p_{\Ai}(t-s, y|x)$ with (\ref{eqn:pAi}), 
that is
\begin{eqnarray}
\mbK_{\Ai}(s,x;t,y)
&=& 
\sum_{a \in \Ai^{-1}(0)} 
\int_{\sqrt{-1} \, \R} \frac{dz}{\sqrt{-1}} \,
p_{\Ai}(s,x|a)
\frac{1}{z-a} \frac{\Ai(z)}{\Ai'(a)}
p_{\Ai}(-t,z|y)
\nonumber\\
&& \quad - {\bf 1}(s>t)p_{\Ai}(s-t,x|y).
\label{eqn:K3a}
\end{eqnarray}

\noindent {\rm (ii)} 
Let $(\Xi_{\cA}(t), \bP_{\Ai})$ be the process, 
which is determinantal with the extended Airy kernel
(\ref{eqn:AiryK}). Then
\begin{equation}
(\Xi_{\cA}(t+\theta), \P_{\Ai}) \to (\Xi_{\cA}(t), \bP_{\Ai})
\quad \mbox{as} \quad \theta\to\infty
\label{relax}
\end{equation}
weakly in the sense of finite dimensional distributions.
\end{thm}
\vskip 3mm

\subsection{Remarks on Extensions of the Results}

(1) \,
By definition (\ref{eqn:MAN}), $M_{\cA}(\xi_{\cA})=0$.
The asymptotic property of the zeros (\ref{eqn:ainf})
implies
\begin{equation}
\zeta^{\cA}(\alpha) \equiv 
\Big( M_{\alpha}(\xi_{\cA}) \Big)^{\alpha}
=\sum_{a \in \cA} \frac{1}{|a|^{\alpha}} < \infty,
\quad \mbox{if} \quad 
\alpha > \widehat{\rho}_f=\frac{3}{2}.
\label{eqn:Azeta}
\end{equation}
In general, order of growth $\rho_f$ of 
a canonical produce (\ref{eqn:Pip}) is
equal to the {\it convergence exponent} $\widehat{\rho}_f$
of the sequence of its zeros \cite{Lev96}.
For $\Ai(z)$, $\rho_f=3/2$.
The function $\zeta^{\cA}(\alpha)$ may be called
the {\it Airy zeta function} \cite{VS04},
which is meromorphic in the whole of $\C$ \cite{FL01}.
Moreover, we know
\begin{equation}
\zeta^{\cA}(2)=M_1(\xi_{\cA}^{\langle 2 \rangle})
=\sum_{a \in {\cA}} \frac{1}{a^2}
=d_1^2 < \infty
\label{eqn:M1A}
\end{equation}
with (\ref{eqn:d0d1}).
Then $\xi_{\cA}$ satisfies the conditions 
({\bf C.1}) and ({\bf C.2}) : 
$\xi_{\cA} \in \mX^{\cA}$.
Since $\xi_{\cA} \in \mM_0$,
Theorem \ref{Theorem:Infinite} guarantees the well-definedness
of the infinite particle system $(\Xi_{\cA}(t), \P_{\cA}^{\xi_{\cA}})$.
(Its equivalence with $(\Xi(t), \P_{\Ai})$ 
is stated in Theorem \ref{Theorem:from_roots} (i).)
Note that the negative divergence (\ref{eqn:Dinf})
of the drift term $D_{\cA_N} t$ of (\ref{eqn:finite_Airy})
in $N \to \infty$ for $t < \infty$ corresponds to that
$\zeta^{\cA}(1)=-\sum_{a \in \cA}(1/a)=\infty$.
This fact and (\ref{eqn:M1A}) mean that
the Airy function has {\it genus} 1 \cite{Lev96,VS04}.

Examples of infinite particle configurations
in $\mX^{\cA}_0$ other than $\xi_{\cA}$ 
are given as follows.
For $\kappa>0$, we put
$$
g^\kappa(x) ={\rm sgn}(x) |x|^{\kappa}, \ x\in\R,
\mbox{ and } 
\eta^{\kappa}(\cdot)=\sum_{\ell\in\Z} 
\delta_{g^\kappa(\ell)}(\cdot).
$$
For any $\kappa > 1/2$ we can confirm by simple calculation
that any configuration $\xi \in \mM_0$ with
$\supp \xi \subset \supp \eta^{\kappa}
=\{g^{\kappa}(\ell):  \ell \in \Z\}$ satisfies ({\bf C.2})(i) 
with any $\alpha \in (1/\kappa, 2)$ 
and some $C_1=C_1(\alpha)>0$ depending on $\alpha$
and ({\bf C.2})(ii) 
with any $\beta \in (0, 2 \kappa-1)$ 
and some $C_2=C_2(\beta) >0$ depending on $\beta$.
Assume that $\xi \in \mM_0$ is chosen so that
$\supp \xi \subset \supp \eta^{\kappa}$
for some $\kappa > 1/2$
and $|M_{\cA}(\xi)| < \infty$.
Then $\xi \in \mX^{\cA}_0$.
The fact (\ref{eqn:ainf}) implies that
this assumption can be satisfied only if
$\kappa \in (1/2, 2/3]$.

\vskip 0.3cm
\noindent
(2) \,
If there exists, however, $\beta' < (\beta-1) \wedge (\beta/2)$ 
for $\xi \in \mM_0$ such that
$\sharp \{ x \in \xi: 
\xi([x-|x|^{\beta'}, x+|x|^{\beta'}]) \geq 2\} = \infty$, 
then $\xi$ does not satisfy the condition
({\bf C.2}) (ii).
In order to include such initial configurations as well as
those with multiple points
in our study of the process $(\Xi_{\cA}(t), \P_{\cA}^{\xi})$
with $\xi(\R)=\infty$,
we introduce another condition for configurations:

\vskip 3mm

\noindent ({\bf C.3})
there exists $\kappa \in (1/2,2/3]$ and $m\in\N$ such that
$$
m(\xi,\kappa)\equiv 
\max_{k\in\Z} \xi\bigg( [ g^\kappa(k), g^\kappa(k+1)] \bigg) \le m.
$$

\vskip 3mm

\noindent We denote by $\mY^{\cA}_{\kappa, m}$ 
the set of configurations $\xi$
satisfying ({\bf C.1}) and ({\bf C.3}) 
with $\kappa\in (1/2,2/3]$ and $m\in\N$, and put
$$
\mY^{\cA} = \bigcup_{\kappa\in (1/2,2/3]}\bigcup_{m\in\N}
\mY^{\cA}_{\kappa, m}.
$$
Noting that the set $\{\xi \in \mM: m(\xi,\kappa)\le m\}$ is 
relatively compact 
for each $\kappa\in (1/2,2/3]$ and $m\in \N$,
we see that $\mY^{\cA}$ is locally compact.

In the present paper, we report our study of
the relaxation process $(\Xi(t), \P_{\Ai})$
from a special initial configuration $\xi_{\cA}$
to the stationary state $\mu_{\Ai}$.
We expect that $\mu_{\Ai}$ is an attractor 
in the configuration space $\mY^{\cA}$ and 
$\xi_{\cA}$ is a point included in the basin.
Motivated by such consideration, 
we are interested in the continuity of the process
with respect to initial configuration.
We have found, however, that
if $\xi(\R)=\infty$, the weak convergence of processes
in the sense of finite dimensional distributions 
can not be concluded from the convergence of
initial configurations in the vague topology.
Following the idea given by our previous paper \cite{KT08},
we introduce a stronger
topology for $\mY^{\cA}$.

Suppose that $\xi, \xi_n \in \mY^{\cA}, n \in \N$.
We say that $\xi_n$ converges 
$\Phi_{\cA}$-{\it moderately} to $\xi$, if
\begin{equation}
\lim_{n\to\infty} \Phi_{\cA}(\xi_n,\sqrt{-1},\cdot)
= \Phi_{\cA}(\xi,\sqrt{-1},\cdot)
\mbox{ uniformly on any compact set of $\C$.}
\label{eqn:THC2}
\end{equation}
It is easy to see that (\ref{eqn:THC2}) is satisfied, 
if the following two conditions hold:
\begin{eqnarray}
&&\lim_{L\to\infty}\sup_{n>0 }\lim_{M\to\infty}
\Bigg| 
M_{\cA}(\xi_n, M)- M_{\cA}(\xi_n, L)
\Bigg|=0,
\label{eqn:THC3}
\\
&&\lim_{L\to\infty} \sup_{n>0 } 
\Bigg| M_2(\xi_n)- M_2(\xi_n, L)
\Bigg|=0.
\label{eqn:THC4}
\end{eqnarray}
By the similar argument given in \cite{KT08},
the following statements are proved. \\
{\rm (i)} If $\xi\in\mY^{\cA}$, the process
$(\Xi_{\cA}, \P_{\cA}^{\xi})$
is well defined. 

\noindent {\rm (ii)} 
Suppose that $\xi, \xi_n \in \mY^{\cA}_{\kappa, n}, n \in\N$,
for some $\kappa\in (1/2,2/3]$ and $m\in\N$.
If $\xi_n$ converges $\Phi_{\cA}$-moderately to $\xi$, then
$
(\Xi_{\cA}, \P_{\cA}^{\xi_n})
\to
(\Xi_{\cA}, \P_{\cA}^{\xi})$
weakly in the sense of finite dimensional distributions
as $n \to \infty$ 
in the vague topology.

Moreover, we can show $\mu_{\Ai}(\mY^{\cA})=1$.
By this fact and the above mentioned 
continuity with respect to initial configurations, 
we can prove that the stationary process
$(\Xi_{\cA}(t), \bP_{\Ai})$, 
which is determinantal
with the extended Airy kernel (\ref{eqn:AiryK}),
is {\it Markovian} \cite{KT09}.

\vskip 0.3cm
\noindent
(3) \,
As mentioned in Introduction, the purpose of
the present paper is to give a method 
for {\it asymmetric} initial configurations
to construct infinite particle systems of
Brownian motions interacting through pair force $1/x$.
In order to clarify the results, we have concentrated
on the case in this paper such that
the initial configuration is $\xi_{\cA}$
(Theorem \ref{Theorem:from_roots})
or its modification $\xi \in \mX^{\cA}_0$; 
see the condition {\bf (C.1)} with (\ref{eqn:MAN})
(Theorem \ref{Theorem:Infinite}).
In the former case the constructed infinite particle
system $(\Xi_{\cA}(t), \P_{\Ai})$ has the
stationary measure $\mu_{\Ai}$, which is 
obtained in the soft-edge scaling limit of the
eigenvalue distribution in GUE well-studied
in the random matrix theory.
Thus we have specified the entire function used
in our analysis in the form $\Phi_{\cA}(\xi,z)$
given by (\ref{eqn:entire_A0}), which is
suitable for the Airy function 
(see (\ref{eqn:relation1})).
The point of our method is to put the relationship
between the entire function appearing in the
correlation kernel $\mbK^{\xi}$, 
the ``typical" initial configuration $\xi_{f}$,
and the drift term in the SDEs providing
finite-particle approximations.
By the same argument as reported here,
the following will be proved.
Let $f$ be the entire function such that $f(0) \not=0$,
it is expressed by the Weierstrass canonical product
of genus one, $\Pi_1(\xi_f,z)$,  
and the zeros can be labelled as 
$0 < |x_1| < |x_2| < \cdots$.
Then with $\xi_f=\sum_{j=1}^{\infty} \delta_{x_j}$
we put
$$
D(\xi_f,N)=\sum_{j=1}^N \frac{1}{x_j}
$$
and introduce the $N$-particle system
$$
Y^f_j(t)=X_j(t)+D(\xi_f, N) t, \quad
1 \leq j \leq N, \quad t \in [0, \infty),
$$
where $\X(t)=(X_1(t), \dots, X_N(t))$
is the solution of Dyson's model (\ref{eqn:Dyson})
starting from the first $N$ zeros of $f$,
$X_j(0)=x_j, 1 \leq j \leq N$.
Then $\Xi_f(t)=\sum_{j=1}^{N} \delta_{Y^f_j(t)}$
converges to the dynamics in $N \to \infty$,
in the sense of finite dimensional distribution,
which is determinantal with the correlation kernel
\begin{eqnarray}
\mbK_f^{\xi_f}(s,x ;t,y)
&=&\int_{\R}\xi_f(dx') 
\int_{\sqrt{-1} \, \R} 
\frac{dy'}{\sqrt{-1}} 
p_{\sin}(s, x-x')
\Phi_1(\xi_f,x',y')
p_{\sin}(-t, y'-y)
\nonumber\\
&& - {\bf 1}(s>t)p_{\sin}(s-t, x-y).
\end{eqnarray}
Moreover, even if the initial configuration $\xi$
is different from $\xi_f$, but it satisfies
the condition
$$
\left| \int_{-L}^{L} \frac{1}{x}
(\xi_f-\xi)(dx) \right| < C_0
$$
for any $L > 0$ with a positive finite
$C_0$ independent of $L$, then the process starting from $\xi$,
is well-defined.
In general, the obtained dynamics with an infinite number
of particles is not stationary, while
Theorem \ref{Theorem:from_roots} gave
the example which converges to a stationary dynamics
$(\Xi_{\cA}(t), \bP_{\Ai})$ in the long-term limit.

\SSC{Properties of the Airy Functions}\label{chap: Airy}

\subsection{Integrals}
By the fact $\Ai''(x)=x \Ai(x)$, 
the following primitive is obtained for $c \not=0$ \cite{VS04},
\begin{eqnarray}
\label{eqn:prim1}
&&\int du \, ( \Ai(c(u+x)) )^2
=(u+x) ( \Ai(c(u+x)) )^2
-\frac{1}{c} ( \Ai'(c(u+x)) )^2, \\
&&\int du \, \Ai(c(u+x)) \Ai(c(u+y))
\nonumber\\
\label{eqn:prim2}
&& \qquad 
= \frac{\Ai'(c(u+x)) \Ai(c(u+y))
- \Ai(c(u+x)) \Ai'(c(u+y))}
{c^2(x-y)}.
\end{eqnarray}
By setting $c=1$ and integral interval be
$[0, \infty)$ in (\ref{eqn:prim2}), we obtain
the integral
$$
\int_{0}^{\infty} du \, \Ai(u+x) \Ai(u+y)
=\frac{\Ai(x) \Ai'(y)-\Ai'(x) \Ai(y)}
{x-y},
$$
since $\lim_{x \to \infty} \Ai(x)=\lim_{x \to \infty} \Ai'(x)=0$
by (\ref{eqn:Aiasym}). 
If we set $y = a \in {\cA}$ and $x=z \not=a$, then
$$
\int_{0}^{\infty} du \, 
\Ai(u+z) \Ai(u+a)
=\frac{\Ai(z) \Ai'(a)}{z-a},
$$
since $\Ai(a)=0$.
Then we have the expression
\begin{equation}
\frac{1}{z-a}
\frac{\Ai(z)}{\Ai'(a)}
=\frac{1}{(\Ai'(a))^2}
\int_{0}^{\infty} du \, 
\Ai(u+z) \Ai(u+a)
\label{eqn:exp1}
\end{equation}
for $a \in {\cA}, z \not=a$.

\subsection{Airy transform}

The following integral formulas are proved 
\cite{Joh03}.

\begin{lem}
\label{thm:eAiAi}
For $c > 0, x, y \in \R$
\begin{eqnarray}
\label{eqn:eAiAi}
&& \int_{\R} du \, e^{c u} \Ai(u+x) \Ai(u+y)
=\frac{1}{\sqrt{4 \pi c}}
e^{-(x-y)^2/(4c)-c(x+y)/2+c^3/12}, \\
\label{eqn:eAiAi2}
&& \int_{\R} d y \int_{\R} du \,
e^{c u} \Ai(u+x) \Ai(u+y)
= e^{-cx+c^3/3}.
\end{eqnarray}
\end{lem}
\noindent{\it Proof.} 
Consider the integral
\begin{eqnarray}
I &=& \int_{\R} du \, e^{c x}
\Ai(u+x) \Ai(u+y) \nonumber\\
&=& \int_{0}^{\infty} du \, e^{c u}
\Ai(u+x) \Ai(u+y)
+ \int_{-\infty}^{0} du \, e^{c u}
\Ai(u+x) \Ai(u+y).
\nonumber
\end{eqnarray}
By the definition of the Airy function (\ref{eqn:Airy1}),
for any $\eta >0$,
\begin{eqnarray}
I &=& \int_{0}^{\infty} du \, e^{c u}
\frac{1}{(2 \pi)^2} \int_{\Im z = \eta} dz \,
\int_{\Im w > c - \eta} d w \,
e^{\sqrt{-1} \{ z^3/3+(u+x)z +w^3/3+(y+u)w \}}
\nonumber\\
&+& \int_{-\infty}^{0} du \, e^{c u}
\frac{1}{(2 \pi)^2} \int_{\Im z = \eta} dz \,
\int_{\Im w < c - \eta} d w
e^{\sqrt{-1}\{ z^3/3+ (u+x)z+w^3/3+(u+y)w \}}
\nonumber\\
&=& \frac{1}{(2\pi)^2} \int_{\Im z = \eta} dz 
\int_{\Im w > c - \eta} dw \,
e^{\sqrt{-1}(z^3/3+ x z + w^3/3+ y w)}
\int_{0}^{\infty} du \,
e^{ \{c+\sqrt{-1}(z+w)\} u } \nonumber\\
&+& \frac{1}{(2\pi)^2} \int_{\Im z = \eta} dz 
\int_{\Im w < c - \eta} dw \,
e^{\sqrt{-1}(z^3/3+x z+ w^3/3+yw)}
\int_{-\infty}^{0} du \,
e^{ \{ c+\sqrt{-1}(z+w)\} u } \nonumber\\
&=& - \frac{1}{(2\pi)^2} \int_{\Im z = \eta} dz 
\int_{\Im w > c - \eta} dw \,
e^{\sqrt{-1}(z^3/3+xz + w^3/3+ yw)}
\frac{1}{c+\sqrt{-1}(z+w)} \nonumber\\
&+& \frac{1}{(2\pi)^2} \int_{\Im z = \eta} dz 
\int_{\Im w < c - \eta} dw \,
e^{\sqrt{-1}(z^3/3+xz + w^3/3+yw)}
\frac{1}{c+\sqrt{-1}(z+w)}
\nonumber\\
&=& \frac{1}{(2 \pi)^2} \int_{\Im z =\eta} dz \,
e^{\sqrt{-1}(z^3/3+xz)} 
\left\{ \int_{\Im w < c-\eta} d w 
- \int_{\Im w > c-\eta} dw \right\}
\frac{e^{\sqrt{-1}(w^3/3+yw)}}{c+\sqrt{-1}(z+w)}.
\nonumber
\end{eqnarray}
It is equal to the integral
$$
\frac{1}{2\pi} \int_{\Im z = \eta} dz \,
e^{\sqrt{-1}(z^3/3+xz)} 
\frac{1}{2 \pi \sqrt{-1}} 
\oint_{C} dw \, \frac{e^{\sqrt{-1}(w^3/3+\sqrt{-1}yw)}}
{w-(\sqrt{-1} c-z)},
$$
where $C$ is a closed contour on $\C$
encircling a pole at $w=\sqrt{-1} c-z$
once in the positive direction.
By performing the Cauchy integral, 
we have
\begin{eqnarray}
I &=& \frac{1}{2 \pi} \int_{\Im z = \eta} dz \,
e^{\sqrt{-1}(z^3/3+xz)}
e^{\sqrt{-1} (\sqrt{-1} c-z)^3/3
-y (c+\sqrt{-1}z)}
\nonumber\\
&=& \frac{1}{2\pi} e^{c^3/3-c y}
\int_{\Im z = \eta} dz \,
e^{ - c z^2 + \sqrt{-1} (x-y+c^2) z }.
\nonumber
\end{eqnarray}
By performing a Gaussian integral, we have (\ref{eqn:eAiAi}).
Since
$$
-\frac{1}{4 c}(x-y)^2 - \frac{c}{2}(x+y)
+\frac{c^3}{12}
= - \frac{1}{4 c} \Big\{ y -(x-c^2) \Big\}^2
- c x  + \frac{c^3}{3},
$$
the Gaussian integration of (\ref{eqn:eAiAi})
with respect to $y$ gives (\ref{eqn:eAiAi2}). \qed
\vskip 0.5cm

By setting $c=t/2 >0$ in (\ref{eqn:eAiAi2})
and $c=(t-s)/2 >0$ in (\ref{eqn:eAiAi}), respectively, 
we obtain the equalities
\begin{eqnarray}
&& \int_{\R} dy \int_{\R} du \, e^{ut}
\Ai(u+x) \Ai(u+y)= g(t,x),
\nonumber\\
&& \int_{\R} du \,
e^{u(t-s)/2} \Ai(u+x) \Ai(u+y) 
= \frac{g(t,y)}{g(s,x)} q(s,t, y-x)
\nonumber
\end{eqnarray}
with (\ref{eqn:g1}) and (\ref{eqn:q}).
Thus we have defined $g(t,x)$ for any $ t \in \R$
by (\ref{eqn:g1}) and 
$p_{\Ai}(s,x;t,y)$ for any $s, t \in \R, s \not= t$ by (\ref{eqn:pAi}).
As special cases, we have
\begin{eqnarray}
\label{eqn:equalM1}
p_{\Ai}(t,y|x)
&=& g(t,y) q(0,t,y-x), \\
\label{eqn:equalM2}
p_{\Ai}(-t,y|x) &=& \frac{1}{g(t,x)} q(t,0,y-x),
\quad t > 0, x,y, \in \R.
\end{eqnarray}

If we take the $c \to 0$ limit 
in (\ref{eqn:eAiAi2}), we obtain
$$
\int_{\R} d \xi \,
\int_{\R} d x \,
\Ai(\xi-x) \Ai(\xi'-x) = 1.
$$
The expression (\ref{eqn:eAiAi}) and the above
result implies
the orthonormality of the Airy function in the sense;
\begin{equation}
\left(\int_{\R} du \,
\Ai(u+x) \Ai(u+y) \right) dy =\delta_x(dy).
\label{eqn:ortho}
\end{equation}
The {\it Airy transform} $f(x) \mapsto \varphi(\xi)$
is then defined by
\begin{equation}
\varphi(\xi)=\int_{\R} dx \,
f(x) \Ai(\xi+x),
\label{eqn:Atran}
\end{equation}
and the inverse transform is given by
$\displaystyle{
f(x)=\int_{\R} d \xi \, \varphi(\xi)
\Ai(\xi+x).
}$
Now a parameter $c \in \C$
is introduced and the family of functions are defined as
$\{
w_{c}(x)= \Ai(x/c)/|c|
\}$. 
The Airy transform (\ref{eqn:Atran}) is then
generalized as
$$
\varphi_{c}(\xi) =
\int_{\R} d \xi \, f(x) w_{c}(\xi+x)
= \frac{1}{|c|} \int_{\R} dx f(x)
\Ai \left( \frac{\xi+x}{c} \right).
$$
\begin{lem}
\label{thm:AiryTrans}
The Airy transform with $c$ of the
normalized Gaussian function
$
  f(x)=e^{-x^2}/\sqrt{\pi}
$
is given by
$\varphi_{c}(\xi)=|c|^{-1}
e^{\{\xi+1/(24 c^3)\}/(4c^3)}
\Ai(\xi/c+1/(16 c^4))$.
That is, 
\begin{equation}
 \int_{\R} dx \, 
\frac{1}{\sqrt{\pi}} e^{-x^2}
\Ai \left( \frac{\xi+x}{c} \right)
= \exp \left\{ \frac{1}{4 c^3}
\left( \xi + \frac{1}{24 c^3} \right) \right\}
\Ai \left( \frac{\xi}{c}
+ \frac{1}{16 c^4} \right).
\label{eqn:AiTrans}
\end{equation}
\end{lem}
\noindent{\it Proof.} 
By the definition of the Airy function (\ref{eqn:Airy1}),
\begin{eqnarray}
I &=& \int_{\R} dx \,
\frac{1}{\sqrt{\pi}} e^{-x^2}
\Ai \left( \frac{\xi+x}{c} \right)
\nonumber\\
&=& \int_{\R} dx \,
\frac{1}{\sqrt{\pi}} e^{-x^2}
\frac{1}{2 \pi} \int_{\R} dk \,
e^{\sqrt{-1}\{k^3/3+(\xi+x)k/c\}}
\nonumber\\
&=& \frac{1}{2\pi} \int_{\R} dk \,
e^{\sqrt{-1} k^3/3+ \sqrt{-1} \xi k/c}
\frac{1}{\sqrt{\pi}} \int_{\R} dx \,
\exp \left( - x^2 - \sqrt{-1} \frac{k}{c} x \right).
\nonumber
\end{eqnarray}
By performing the Gaussian integral we have
$$
I=\frac{1}{2 \pi } \int_{\R} dk \,
\exp \left( \sqrt{-1} \frac{k^3}{3}
+ \sqrt{-1} \frac{\xi k}{c} - \frac{k^2}{4 c^2} \right).
$$
By completing a cube, we find the equality
\begin{eqnarray}
&& \sqrt{-1} \frac{k^3}{3}
+ \sqrt{-1} \frac{\xi k}{c} - \frac{k^2}{4 c^2}
=\sqrt{-1} \frac{1}{3} \left( 
k+\sqrt{-1} \frac{1}{4 c^2} \right)^3
\nonumber\\
&& \qquad
+ \sqrt{-1} \left( \frac{\xi}{c}+\frac{1}{16 c^4} \right)
\left(k + \sqrt{-1} \frac{1}{4 c^2} \right)
+ \frac{1}{4 c^3} 
\left( \xi + \frac{1}{24 c^3} \right).
\nonumber
\end{eqnarray}
By using the definition of the Airy function (\ref{eqn:Airy1}), 
(\ref{eqn:AiTrans}) is obtained.
\qed
\vskip 0.3cm

For $t>0, y, u \in \R$
we will obtain the equality
$$
\int_{\sqrt{-1} \, \R} \frac{dz}{\sqrt{-1}} \,
\Ai(u+z) q(t, 0, z-y) 
=\int_{\R} dx \, \frac{1}{\sqrt{\pi}}
e^{-x^2} \Ai \left( \sqrt{-2t} \, x +
u+y-\frac{t^2}{4} \right)
$$
by changing the integral variable as
$z \mapsto x=(z-y+t^2/4)/\sqrt{-2t}$.
If we set 
$\xi=(u+y-t^2/4)/\sqrt{-2t}$
and $c=1/\sqrt{-2t}$,
the RHS is identified with the LHS of
(\ref{eqn:AiTrans}).
Since
$$
\frac{\xi}{c}+\frac{1}{16c^4}=u+y, \quad
\frac{1}{4 c^3} \left( \xi+\frac{1}{24 c^3} \right)
=\left( -\frac{ty}{2}+\frac{t^3}{24} \right) - \frac{ut}{2},
$$
(\ref{eqn:AiTrans}) of Lemma \ref{thm:AiryTrans}
with (\ref{eqn:g1}) gives 
$$
 \int_{\sqrt{-1} \, \R} \frac{dz}{\sqrt{-1}} \,
\Ai(u+z) q(t, 0, z-y) 
= g(t, y) e^{-ut/2} \Ai(u+y),
\quad t >0, y, u \in \R.
$$
Combination with (\ref{eqn:equalM2}) gives
\begin{equation}
\int_{\sqrt{-1} \, \R} \frac{dz}{\sqrt{-1}} \,
\Ai(u+z) p_{\Ai}(-t, z|y)
= e^{-ut/2} \Ai(u+y),
\quad t >0, y, u \in \R
\label{eqn:equalM3}
\end{equation}

\subsection{Fourier-Airy series}

Let us consider the integral
$$
I_{\ell \ell'} = \int_{0}^{\infty} dx \,
\Ai(x+a_{\ell}) \Ai(x+a_{\ell'}),
\quad a_{\ell}, a_{\ell'} \in {\cA}.
$$
In the case $\ell \not= \ell'$, the formula (\ref{eqn:prim2})
gives
$$
I_{\ell \ell'}
=\frac{\Ai'(a_{\ell}) \Ai(a_{\ell'})
-\Ai(a_{\ell}) \Ai'(a_{\ell'})}{a_{\ell}-a_{\ell'}}=0,
$$
whereas if $\ell =\ell'$, the formula (\ref{eqn:prim1}) gives
$
I_{\ell \ell}= ( \Ai'(a_{\ell}) )^2.
$
Therefore the functions
\begin{equation}
\left\{ \frac{\Ai(x+a_{\ell})}{\Ai'(a_{\ell})},
\ell \in \N \right\}
\label{eqn:Aiset}
\end{equation}
form an orthogonal basis for 
$f \in L^2(0, \infty)$ 
(see Sect. 4.12 in \cite{Tit62}).
The completeness of (\ref{eqn:Aiset}) is also
established :
\begin{equation}
\sum_{\ell \in \N}
\frac{\Ai(x+a_{\ell}) \Ai(y+a_{\ell})}
{( \Ai'(a_{\ell}) )^2} dy
=\delta_{x} (dy),
\quad x,y \in (0,\infty).
\label{eqn:complete}
\end{equation}
Then for any $f \in L^2(0, \infty)$, 
we can write the expression
$$
f(x)=\sum_{\ell \in \N} c_{\ell}
\frac{\Ai(x+a_{\ell})}{\Ai'(a_{\ell})},
\quad x\in [0,\infty),
$$
and call it the {\it Fourier-Airy series expansion}.
The coefficients $c_{\ell}$ of this expansion
are determined by
$
c_{\ell}=\{ 
\int_{0}^{\infty} dx \, f(x) \Ai(x+a_{\ell})\}
/\Ai'(a_{\ell})
$.

\SSC{Proof of Results}

\subsection{Proof of Proposition \ref{Proposition:Finite}}
With (\ref{eqn:DA}) we put
$$
\widehat{g}(s,x)=\exp\left\{
-D_{\cA_N} \left(
\frac{D_{\cA_N} s}{2}+\frac{s^2}{4}-x \right)\right\},
\, s,x\in\R.
$$
By the definition (\ref{eqn:q}) of $q$,
for $s >0, x, x' \in \R$, we have
\begin{eqnarray}
&&\frac{p_{\sin}(s,(x-D_{\cA_N}s-s^2/4)-x')}
{q(0, s, x-x')}
\nonumber\\
&& \qquad =\exp\left[
-\frac{1}{2s}
\left\{
\left(x-D_{\cA_N}s -\frac{s^2}{4}-x' \right)^2
-\left(x-\frac{s^2}{4}-x' \right)^2 \right\} \right]
\nonumber\\
&& \qquad =\exp\left[
-D_{\cA_N} s \left( 
\frac{D_{\cA_N}s}{2} + \frac{s^2}{4}-x+x' \right) \right]
\nonumber\\
&& \qquad =\widehat{g}(s,x) e^{-D_{\cA_N} x'},
\nonumber
\end{eqnarray}
and for $t > 0, y, y' \in \R$, we have
\begin{equation}
\frac{p_{\sin}(-t, y'-(y-D_{\cA_N}t-t^2/4))}
{q(t,0, y'-y)}
=\frac{1}{\widehat{g}(t,y)} e^{D_{\cA_N} y'}.
\nonumber
\end{equation}
Similarly, we have
\begin{eqnarray}
&&\frac{p_{\sin}(s-t,(x-D_{\cA_N} s-s^2/4)-
(y-D_{\cA_N} t-t^2/4))}{q(t, s, y-x)}
\nonumber\\
&&=\exp\left[
-\frac{1}{2(s-t)}
\left\{
\left(x-y-D_{\cA_N}(s-t)-\frac{s^2-t^2}{4} \right)^2
\right. \right. \nonumber\\
&& \qquad \qquad \left. \left.
-\left( \left(x-\frac{s^2}{4} \right)
-\left(y-\frac{t^2}{4} \right) \right)^2
\right\}
\right]
\nonumber\\
&&=\exp\left[- D_{\cA_N}
\left\{ \frac{D_{cA_N}}{2}(s-t)
+\frac{s^2-t^2}{4}-(x-y) \right\} \right]
=\frac{\widehat{g}(s,x)}{\widehat{g}(t,y)}.
\nonumber
\end{eqnarray}
Then we have
\begin{eqnarray}
&& \mbK^{\xi^N}\left(s, x-D_{\cA_N}s-\frac{s^2}{4}; 
t, y-D_{\cA_N}t- \frac{t^2}{4} \right)
\nonumber\\
&&
= \frac{\widehat{g}(s,x)}{\widehat{g}(t,y)}
\left[ 
\int_{\R} \xi^{N}(dx') 
\int_{\sqrt{-1} \, \R} 
\frac{dy'}{\sqrt{-1}} \, \right.
\nonumber\\
&& \qquad \qquad \times
q(0, s, x-x') 
e^{-D_{\cA_N}x'} \Phi_0(\xi^N, x', y') e^{D_{\cA_N}y'} 
q(t, 0, y'-y) \nonumber\\
&& \hskip 6cm
-{\bf 1}(s > t) q(t, s, x-y) \Bigg].
\nonumber
\end{eqnarray}
The identity (\ref{eqn:entire_A}) implies
$$
e^{-D_{\cA_N}x'} \Phi_0(\xi^N,x',y')
e^{D_{\cA_N}y'}
=\Phi_{\cA}(\xi^N,x',y').
$$
By the gauge invariance of determinantal processes
(Lemma \ref{thm:gauge}), the proof is completed.
\qed

\subsection{Proof of Theorem \ref{Theorem:Infinite}}

First we prepare some lemmas for proving 
Theorem \ref{Theorem:Infinite}.
\begin{lem}\label{thm:4_1}
Let $\alpha\in (1,2)$ and $\delta> \alpha-1$.
Suppose that $M_\alpha(\xi)<\infty$ and put 
$L_0=L_0(\alpha, \delta, \xi)= 
(2M_\alpha(\xi))^{\alpha/(\delta-\alpha+1)}$.
Then
$$
M_1 (\xi,L) \le L^\delta, \quad L\ge L_0.
$$
\end{lem}
Since this lemma was proved as Lemma 4.3 in 
\cite{KT08}, here we omit the proof.

\begin{lem}\label{thm:4_2}
If $\xi$ satisfies  {\rm ({\bf C.2})} {\rm (i)} and {\rm (ii)},
for any $\theta\in (\alpha \vee (2-\beta), 2)$
there exists $C=C(C_1,C_2,\theta)>0$ such that
\begin{equation}
\left| 
\int_{\{0,a\}^{\rm c}} \
\left( \frac{1}{x} -\frac{1}{x-a}\right)
\xi(dx)\right|
\le  C|a \vee 1|^{\theta-1},
\quad a\in \supp(\xi-\delta_0).
\label{Difference}
\end{equation}
\end{lem}
\vskip 3mm
\noindent {\it Proof.}
We divide $\{0,a\}^{\rm c}$ into three sets 
$A_1= \{x\in\{0,a\}^{\rm c} : |x|\le |a|/2 \}$,
$A_2= \{x\in\{0,a\}^{\rm c} : |a|/2 <|x|\le 2|a| \}$,
and
$A_3= \{x\in\{0,a\}^{\rm c} : 2|a|<|x|\}$, 
and put
$$
I_j= \int_{A_j}\left| \frac{1}{x} -\frac{1}{x-a}
\right|\xi(dx)
\quad j=1,2,3.
$$
When $x\in A_1$, $|x^2-a^2|\ge 3 a^2/4$ 
and $|x+a|\le 3 |a|/2$,
and then
$$
I_1 = \int_{A_1}\frac{|a||x+a|}{|x||x^2-a^2|}\xi(dx)
\le 2M_1\left(\xi,\frac{|a|}{2}\right).
$$
By Lemma \ref{thm:4_1} for any $\delta>\alpha-1$,
we can take $C>0$ such that
\begin{equation}
I_1\le C|a \vee 1|^\delta.
\label{I_1}
\end{equation}
When $x\in A_2$, $|x+a| \leq 3|a|$ and $|a|/|x|\le 2$, 
and then
$$
I_2 = \int_{A_2}\frac{|a||x+a|}{|x||x^2-a^2|}\xi(dx)
\le 6|a|M_1(\tau_{-a^2}\xi^{\langle 2 \rangle}).
$$
From the condition ({\bf C.2}) {\rm (ii)}
\begin{equation}
I_2\le 6 C_2 |a\vee 1|^{1-\beta}.
\label{I_2}
\end{equation}
When $x\in A_3$, $|x-a|> |x|/2$, and then
$$
I_3 = \int_{A_3}\frac{|a|}{|x||x-a|}\xi(dx)
\le 2^{\alpha-1} |a|^{\alpha-1} M_{\alpha}(\xi)^{\alpha}.
$$
From the condition ({\bf C.2}) {\rm (i)}
\begin{equation}
I_3\le 2^{\alpha-1} C_1 |a\vee 1|^{\alpha-1}.
\label{I_3}
\end{equation}
Combining the estimates (\ref{I_1}), (\ref{I_2}) and (\ref{I_3}),
we have (\ref{Difference}).
\qed

\begin{lem}
\label{thm:4_3}
{\rm (i)}\quad
If $\xi$ satisfies the conditions
{\rm ({\bf C.2})} {\rm (i)} and {\rm (ii)},
for any $\theta\in (\alpha \vee (2-\beta), 2)$
there exists $C=C(C_1,C_2,\theta)>0$ such that
\begin{equation}
|\Phi_1(\xi,a,\sqrt{-1} y)|\le 
\exp\left[ C\{(|y|^{\theta} \vee 1) + (|a|^{\theta} \vee 1) \}\right],
\end{equation}
for $y\in\R$ and $a\in\supp \xi$.

\noindent {\rm (ii)}\quad
If $\xi$ satisfies the conditions
{\rm ({\bf C.1}), ({\bf C.2})} {\rm (i)} and {\rm (ii)},
for any $\theta\in (\alpha \vee (2-\beta), 2)$
there exists $C=C(C_0, C_1,C_2,\theta)>0$ such that
\begin{equation}
|\Phi_{\cA}(\xi,a, \sqrt{-1} y)|\le 
\exp\left[ C\{(|y|^{\theta} \vee 1) 
+ (|a|^{\theta} \vee 1) \}\right],
\end{equation}
for $y\in\R$ and $a\in\supp \xi$.
\end{lem}
\vskip 3mm
\noindent {\it Proof.}
We prove {\rm (i)} of this lemma.
From the condition ({\bf C.1}) and the relation 
(\ref{eqn:entire_A}), {\rm (ii)} is easily derived from
{\rm (i)}.
We first consider the case that $a=0\in\supp \xi$.
Remind that
$$
\Phi_1(\xi, 0, z)
= \Pi_{1}(\xi, z)
= \exp\left[
\int_{\{0\}^{\rm c} } \left\{
\log \left(1-\frac{z}{x}\right) + \frac{z}{x} \right\}
\xi(dx)\right].
$$
When $2|z|<|x|$, by using the expansion
$$
\log \left(1-\frac{z}{x}\right)
= - \sum_{k \in \N} \frac{1}{k}
\left(\frac{z}{x}\right)^k,
$$
we have 
$$
\left|\log \left(1-\frac{z}{x}\right) + \frac{z}{x}\right|
\le \left|\frac{z}{x}\right|^2.
$$
Then
\begin{equation}
|\Pi_1(\xi\cap [-2|z|,2|z|]^{\rm c},z)|
\le \exp\left\{|z|^2 \int_{|x|>2|z|} \frac{1}{x^2}\xi(dx) \right\}
\leq \exp \Big\{|z|^{\alpha} M_{\alpha}(\xi)^\alpha \Big\}.
\label{estimate:Phi_1_1}
\end{equation}
On the other hand
$ |1-z/x| \leq e^{|z|/|x|}$.
Then
\begin{equation}
|\Pi_1(\xi\cap [-2|z|,2|z|],z)|
\le \exp\left\{2|z| \int_{|x|\le 2|z|} \frac{1}{|x|}\xi(dx) \right\}
=\exp \Big\{2|z| M_1(\xi,2|z|) \Big\}.
\label{estimate:Phi_1_2}
\end{equation}
From (\ref{estimate:Phi_1_1}) and (\ref{estimate:Phi_1_2}),
with the condition ({\bf C.2}) {\rm (i)} and
Lemma \ref{thm:4_1},
we see that 
for any $\theta\in (\alpha \vee (2-\beta), 2)$
there exists $C=C(C_1,\theta)>0$ such that
\begin{equation}
|\Phi_1(\xi, 0, z)|\le 
\exp\left[ C\{(|z|^{\theta} \vee 1) \}\right],
\label{estimata:a=0}
\end{equation}
for $z\in\C$ and $a\in\supp \xi$.

Next we consider the case that $a\in\supp \xi$ and $a\not=0$.
By the conditions ({\bf C.1}) and ({\bf C.2})
the equality (\ref{eqn:finite_zeta}) is valid.
By (\ref{estimata:a=0})
$$
|\Pi_1(\xi,z)\Pi_1(\xi\cap \{-a\}^{\rm c},-a)|
\le \exp \Big[ C\{(|z|^{\theta} \vee 1)
 + (|a|^{\theta} \vee 1) \} \Big].
$$
By the condition ({\bf C.2}) {\rm (ii)}
$$
|\Phi_0(\xi^{\langle 2 \rangle} \cap \{0\}^{\rm c}, a^2,0)|
\le \exp \Big\{|a|^2 M_1(\tau_{-a^2}\xi^{\langle 2 \rangle}) \Big\}
\le \exp \Big\{C_2 (|a|\vee 1)^{2-\beta} \Big\},
$$
and $|(\sqrt{-1} y/a)^{\xi(\{0\})} a/(a-\sqrt{-1} y)| \leq 1$.
Now we evaluate $S(\xi,a,z)$.
\begin{eqnarray}
&&S(\xi,a,z)= 
\int_{\{0,a\}^{\rm c}}
\left( \frac{z-a}{x-a}-\frac{z-a}{x}\right) \xi(dx)
\nonumber\\
&& \hskip 4cm
+\frac{(z-a)}{-a}\xi(\{0\})
-\frac{z}{a}+\frac{a}{a}-\frac{a}{-a} \xi(\{-a\})
\nonumber\\
&&=(z-a) \int_{\{0,a\}^{\rm c}}
\left(\frac{1}{x-a}-\frac{1}{x}\right) \xi(dx)
-\frac{(1+\xi(\{0\}))z}{a}+1+\xi(\{0\})+\xi(\{-a\}).
\nonumber
\end{eqnarray}
From Lemma \ref{thm:4_2} and the fact 
$1/a^2 \le C_2$ and then 
$|2z/a|\le 2\sqrt{C_2}|z|$, 
we have
$$
|S(\xi,a,z)|\le C|z-a||a\vee 1|^{\theta-1} +2\sqrt{C_2}|z|+3
\le C' \Big\{(|y|^{\theta} \vee 1) 
+ (|a|^{\theta} \vee 1) \Big\}
$$
for some $C'>0$. This completes the proof.
\qed
\vskip 3mm

\noindent{\it Proof of Theorem \ref{Theorem:Infinite}.} 
Note that $\xi\cap [-L,L]$, $L>0$ and $\xi$ 
satisfy ({\bf C.1}) and ({\bf C.2})
with the same constants $C_0, C_1, C_2$ and indices $\alpha, \beta$.
By virtue of Lemma \ref{thm:4_3} {\rm (ii)}
we see that there exists $C >0$ such that
$$
|\Phi_{\cA} (\xi\cap [-L,L],a, \sqrt{-1} y)|
\le \exp \Big[ C \Big\{(|y|\vee 1)^{\theta} 
+ (|a|\vee 1)^{\theta} \Big\} \Big],
$$
$\forall L>0, \; \forall a\in \supp \xi, \forall y \in \R$.
Since for any $y\in\R$
$$
\Phi_{\cA} (\xi\cap [-L,L],a,\sqrt{-1} y) 
\to \Phi_{\cA} (\xi,a,\sqrt{-1} y), \quad L\to\infty,
$$
we can apply {\it Lebesgue's convergence theorem}
to (\ref{eqn:K3}) and obtain
$$
\lim_{L\to\infty} \mbK_{\cA}^{\xi\cap [-L,L]}
\left(s, x; t, y \right)
=\mbK_{\cA}^{\xi}\left(s, x; t, y\right).
$$
Since for any $(s,t) \in (0, \infty)^{2}$ and 
any finite interval $I \subset \R$
$$
\sup_{x, y \in I} \Big|
\mbK_{\cA}^{\xi \cap [-L, L]}(s, x; t, y) \Big| < \infty,
$$
we can obtain the convergence of generating functions
for multitime correlation functions (\ref{eqn:Gxi});
${\cG}^{\xi \cap [-L, L]}[\chi] \to
{\cG}^{\xi}[\chi]$ as $L \to \infty$.
It implies 
$\P_{\cA}^{\xi \cap [-L, L]} \to 
\P_{\cA}^{\xi}$ as $L \to \infty$
in the sense of finite dimensional distributions.
Then the proof is completed.
\qed

\subsection{Proof of Theorem \ref{Theorem:from_roots} }

(i) It is clear that $\xi_{\cA}$ is an element of $\mX_{\cA}^0$.
(See the item (1) of Sect. 2.4.)
Then by Theorem \ref{Theorem:Infinite} 
$(\Xi_{\cA}(t), \P_{\cA}^{\xi_{\cA}^{N}})
\to (\Xi_{\cA}(t), \P_{\cA}^{\xi_{\cA}})$ as $N \to \infty$
in the sense of finite dimensional distributions,
where $(\Xi_{\cA}(t), \P_{\cA}^{\xi_{\cA}})$
is the determinantal with the correlation kernel
\begin{eqnarray}
\mbK^{\xi_{\cA}}_{\cA}(s,x;t,y)
&=& 
\int_{\R}\xi_{\cA}(da) 
\int_{\sqrt{-1} \, \R} \frac{dz}{\sqrt{-1}}
 \, q(0,s, x-a) 
\Phi_{\cA}(\xi_{\cA}, a, z) q(t,0, z-y)
\nonumber\\
&& \quad - {\bf 1}(s>t)q(t,s, x-y).
\nonumber
\end{eqnarray}
Using the equalities (\ref{eqn:equalM1}) and
(\ref{eqn:equalM2}) 
and the definition (\ref{eqn:pAi}), 
we have 
\begin{eqnarray}
\mbK_{\cA}^{\xi_{\cA}}(s,x;t,y)
&=& \int_{\R} \xi_{\cA}(da) 
\int_{\sqrt{-1} \, \R} \frac{dz}{\sqrt{-1}} \,
\frac{p_{\Ai}(s, x|a)}{g(s,x)}
\frac{1}{z-a} \frac{\Ai(z)}{\Ai'(a)}
g(t,y) p_{\Ai}(-t,z|y)
\nonumber\\
&& \quad 
-{\bf 1}(s>t) \frac{g(t,y)}{g(s,x)} 
p_{\Ai}(s-t,x|y)
\nonumber\\
&=& \frac{g(t,y)}{g(s,x)} 
\mbK_{\Ai}(s,x;t,y)
\nonumber
\end{eqnarray}
with (\ref{eqn:K3a}), 
where we have used (\ref{eqn:relation1b}) of Lemma \ref{thm:2_2}.
By the gauge invariance, Lemma \ref{thm:gauge},
$(\Xi_{\cA}(t), \P_{\cA}^{\xi_{\cA}})
=(\Xi_{\cA}(t), \P_{\Ai})$
in the sense of finite dimensional distributions. \\
(ii) If we use the expression (\ref{eqn:exp1}), 
(\ref{eqn:K3a}) becomes
\begin{eqnarray}
\mbK_{\Ai}(s,x;t,y) &=& 
\sum_{a \in \Ai^{-1}(0)} 
\int_{\sqrt{-1} \, \R} \frac{dz}{\sqrt{-1}} 
\, p_{\Ai}(s, x|a) 
\nonumber\\
&& \qquad \times
\frac{1}{(\Ai'(a))^2} 
\int_{0}^{\infty} du \, \Ai(u+z) \Ai(u+a)
p_{\Ai}(-t, z|y)
\nonumber\\
&& - {\bf 1}(s>t) p_{\Ai}(s-t,x|y).
\nonumber
\end{eqnarray}
By (\ref{eqn:equalM3}), the first term of the RHS equals
$$
\sum_{a \in \Ai^{-1}(0)} p_{\Ai}(s,x|a)
\frac{1}{(\Ai'(a))^2} \int_{0}^{\infty} du \,
e^{-ut/2} \Ai(u+y) \Ai(u+a).
$$
Since $s>0$, we can use the expression (\ref{eqn:pAi0})
for $p_{\Ai}(s, x|a)$ and the above is written as
$$
\int_{0}^{\infty} du \int_{\R} dw \,
e^{-ut/2+ws/2}
\Ai(u+y) \Ai(w+x)
\sum_{\ell \in \N}
\frac{\Ai(u+a_{\ell}) \Ai(w+a_{\ell})}
{( \Ai'(a_{\ell}) )^2}.
$$
From the completeness (\ref{eqn:complete}), the above gives
$$
\mbK_{\Ai}(s,x; t,y)
={\bK}_{\Ai}(t-s, y|x)+R(s,x;t,y)
$$
with the extended Airy kernel $\bK_{\Ai}$ given by
(\ref{eqn:AiryK}) and
\begin{eqnarray}
R(s,x;t,y) &=& \int_{0}^{\infty} du \int_{-\infty}^0 dw \,
e^{-ut/2+ws/2}
\Ai(u+y) \Ai(w+x)
\nonumber\\
&& \quad \times
\sum_{\ell \in \N}
\frac{\Ai(u+a_{\ell}) \Ai(w+a_{\ell})}
{( \Ai'(a_{\ell}) )^2}.
\nonumber
\end{eqnarray}
Since for any fixed $s,t >0$
$\lim_{\theta\to\infty}|R(s+\theta,x;t+\theta,y)|\to 0$
uniformly on any compact subset of $\R^2$, 
(\ref{eqn:convKAi}) holds in the same sense.
Hence we obtain (\ref{relax}).
This completes the proof.
\qed

\vskip 1cm
\begin{small}
\noindent{\it Acknowledgments.} 
M.K. is supported in part by
the Grant-in-Aid for Scientific Research (C)
(No.21540397) of Japan Society for
the Promotion of Science.
H.T. is supported in part by
the Grant-in-Aid for Scientific Research 
(KIBAN-C, No.19540114) of Japan Society for
the Promotion of Science.
\vskip 1cm

\end{small}
\end{document}